\documentclass[a4paper,11pt]{article}
\usepackage[latin1]{inputenc}
\usepackage[francais,english]{babel}    
\usepackage{amssymb}
\usepackage{amsmath}
\usepackage{amsthm}
\usepackage{multicol}
\usepackage{graphicx}
\usepackage{color}
\usepackage[table]{xcolor}
\usepackage[T1]{fontenc}
\usepackage{yfonts}
\usepackage{placeins}
\usepackage{array}
\usepackage{dsfont}
\usepackage{stmaryrd}
\usepackage{verbatim}
\usepackage{caption}
\usepackage{transparent}
\usepackage{bbold}
\usepackage{hyperref}
\usepackage{enumitem}
\usepackage{tikz}
\usepackage{float}
\usetikzlibrary{patterns}
\usetikzlibrary{arrows.meta}
\usetikzlibrary{decorations.markings}
\usetikzlibrary{matrix,arrows,calc,fit,cd,positioning,intersections,arrows.meta,braids}
\usepackage{pgfplots}
\usepackage{tikz-3dplot}

\newtheorem{thm}{Theorem}[section]
\newcommand{\bthm}{\begin{thm}}
\newcommand{\ethm}{\end{thm}}

\newtheorem{thmi}{Theorem}
\newcommand{\bthmi}{\begin{thmi}}
\newcommand{\ethmi}{\end{thmi}}

\newtheorem{cori}[thmi]{Corollary}
\newcommand{\bcori}{\begin{cori}}
\newcommand{\ecori}{\end{cori}}

\newtheorem{mthm}{Theorem}
\newcommand{\bmthm}{\begin{mthm}}
\newcommand{\emthm}{\end{mthm}}

\newtheorem{mcor}[mthm]{Corollary}
\newcommand{\bmcor}{\begin{mcor}}
\newcommand{\emcor}{\end{mcor}}
\newtheorem{mconj}[mthm]{Conjecture}
\newcommand{\bmconj}{\begin{mconj}}
\newcommand{\emconj}{\end{mconj}}
\newtheorem{mpro}[mthm]{Proposition}
\newcommand{\bmpro}{\begin{mpro}}
\newcommand{\empro}{\end{mpro}}

\newtheorem*{conj}{Conjecture}
\newcommand{\bconj}{\begin{conj}}
\newcommand{\econj}{\end{conj}}

\newtheorem*{question}{Question}
\newcommand{\bq}{\begin{question}}
\newcommand{\eq}{\end{question}}

\newtheorem*{thn}{Theorem}
\newcommand{\bthn}{\begin{thn}}
\newcommand{\ethn}{\end{thn}}

\newtheorem{exo}{Exercise}
\newcommand{\bex}{\begin{exo}}
\newcommand{\eex}{\end{exo}}

\newtheorem{sol}{Solution}
\newcommand{\bsol}{\begin{sol}}
\newcommand{\esol}{\end{sol}}

\newtheorem{pro}[thm]{Proposition}
\newcommand{\bpro}{\begin{pro}}
\newcommand{\epro}{\end{pro}}

\newtheorem{cor}[thm]{Corollary}
\newcommand{\bcor}{\begin{cor}}
\newcommand{\ecor}{\end{cor}}

\newtheorem{lem}[thm]{Lemma}
\newcommand{\blem}{\begin{lem}}
\newcommand{\elem}{\end{lem}}

\theoremstyle{definition}

\newtheorem{defi}[thm]{Definition}
\newcommand{\bdf}{\begin{defi}}
\newcommand{\edf}{\end{defi}}

\newtheorem*{defis}{Definition}
\newcommand{\bdfs}{\begin{defis}}
\newcommand{\edfs}{\end{defis}}

\newtheorem*{rmk}{Remark}
\newcommand{\brk}{\begin{rmk} \upshape}
\newcommand{\erk}{\end{rmk}}

\newtheorem*{rmks}{Remarks}
\newcommand{\brks}{\begin{rmks} \upshape}
\newcommand{\erks}{\end{rmks}}

\newtheorem*{exe}{Example}
\newcommand{\bexe}{\begin{exe} \upshape}
\newcommand{\eexe}{\end{exe}}

\newtheorem*{exes}{Examples}
\newcommand{\bexes}{\begin{exes} \upshape}
\newcommand{\eexes}{\end{exes}}

\newtheorem*{pre}{Proof}
\newcommand{\bp}{\begin{pre} \upshape}
\newcommand{\ep}{\hfill \qed \end{pre}}
\newcommand{\epp}{\end{pre}}

\newcommand{\beq}{\begin{eqnarray*}}
\newcommand{\eeq}{\end{eqnarray*}}

\newcommand{\beqn}{\begin{equation}}
\newcommand{\eeqn}{\end{equation}}

\newcommand{\ben}{\begin{enumerate}}
\newcommand{\een}{\end{enumerate}}
\newcommand{\bit}{\begin{itemize} \renewcommand{\labelitemi}{$\bullet$} \renewcommand{\labelitemii}{$\star$}}
\newcommand{\eit}{\end{itemize}}

\newcommand{\bfg}{
\begin{figure}[H]
\begin{center}}
\newcommand{\efg}{
\end{center}
\end{figure}
\FloatBarrier}

\newcolumntype{M}[1]{>{\raggedright}m{#1}}

\newcommand{\R}{\mathbb{R}}
\newcommand{\Q}{\mathbb{Q}}
\newcommand{\N}{\mathbb{N}}
\newcommand{\Z}{\mathbb{Z}}

\newcommand{\K}{\mathbb{K}}

\renewcommand{\H}{\mathbb{H}}

\newcommand{\F}{\mathbb{F}}

\newcommand{\bs}{\symbol{92}}

\newcommand{\ov}{\overline}

\renewcommand{\tilde}{\widetilde}

\renewcommand{\t}{ ^t\!}

\renewcommand{\dim}{\operatorname{dim}}

\newcommand{\Isom}{\operatorname{Isom}}

\renewcommand{\max}{\operatorname{max}}

\newcommand{\Hull}{\operatorname{Hull}}

\newcommand{\diam}{\operatorname{diam}}

\newcommand{\Lip}{\operatorname{Lip}}

\newcommand{\eps}{\varepsilon}
\newcommand{\st}{\, | \,}

\newcommand{\ra}{\rightarrow}

\newcommand{\liml}{\lim\limits}

\newcommand{\f}{\frac}

\renewcommand{\geq}{\geqslant}
\renewcommand{\leq}{\leqslant}

\renewcommand{\log}{\operatorname{log}}

\newcommand{\id}{\operatorname{id}}

\newcommand{\GL}{\operatorname{GL}}
\newcommand{\SL}{\operatorname{SL}}
\newcommand{\SO}{\operatorname{SO}}

\renewcommand{\>}{\rangle}
\newcommand{\pif}{{+\infty}}

\newcommand{\mk}{\medskip}

\makeatletter
\def\Ddots{\mathinner{\mkern1mu\raise\p@
\vbox{\kern7\p@\hbox{.}}\mkern2mu
\raise4\p@\hbox{.}\mkern2mu\raise7\p@\hbox{.}\mkern1mu}}
\makeatother

\makeatletter
\def\maketitles{%
  \null
  \thispagestyle{empty}%
  \vfill
  \begin{center}\leavevmode
    \normalfont
    {\LARGE \@title\par}%
    \vskip 1.2cm
    {\large \@author\par}%
    \vskip 1.2cm
    {\large \@subtitle\par}%
    \vskip 0.8cm
    {\large \@date\par}%
  \end{center}%
  \vfill
  \null
  \cleardoublepage
  }

\def\date#1{\def\@date{#1}}
\def\author#1{\def\@author{#1}}
\def\title#1{\def\@title{#1}}
\def\subtitle#1{\def\@subtitle{#1}}
\makeatother

\title{Group actions on injective spaces and Helly graphs\\ Lecture Notes \\ Minicourse at the CRM Montr\'eal, June 2023}
\author{Haettel, Thomas\footnote{Thomas Haettel, thomas.haettel@umontpellier.fr, IMAG, Univ Montpellier, CNRS, France, and IRL 3457, CRM-CNRS, Universit\'{e} de Montr\'{e}al, Canada.}}

\begin{document}

\maketitle

\tableofcontents

\definecolor{gris1}{gray}{0.75}
\definecolor{gris2}{gray}{0.85}
\rowcolors{2}{gris2}{white}

\section{Introduction}

Injective metric spaces and Helly graphs have recently become a very active object of study in geometric group theory, intiated notably by the brilliant article~\cite{lang} by Lang. Whereas it is a notorious open problem whether Gromov-hyperbolic groups act geometrically on CAT(0) spaces, Lang proves that Gromov-hyperbolic groups act geometrically on injective metric spaces. It turns out that injective metric spaces enjoy many properties which are typical of nonpositive curvature, but are sometimes better behaved than CAT(0) spaces. Indeed it is easier to build injective metric spaces, and one can sometimes deduce stronger properties for groups acting on them.

\mk

For instance, in addition to hyperbolic groups, one can show that plenty of groups with nonpositive curvature flavour have an interesting action on an injective metric space, notably braid groups, mapping class groups and higher rank lattices (see Section~\ref{sec:summary_examples_properties} for a precise list). Moreover, given an isometric action of a group on an injective metric space, one can often deduce strong consequences, such as algorithmic properties, distortion of elements, biautomaticity, nice classifying spaces (see Section~\ref{sec:summary_examples_properties} for many consequences).

\mk

Injective metric spaces and Helly graphs should be thought as reflecting "$L^\infty$ geometry", where CAT(0) spaces reflect "$L^2$ geometry" and metric median spaces reflect $L^1$ geometry. Rough and informal analogies between these three geometries are listed in Section~\ref{sec:L1L2Linfini}.

\mk

The purpose of this article is to survey basic definitions and properties of injective metric spaces and Helly graphs. When available, short proofs are provided. Some results are original. We chose to emphasize the parallel between injective metric spaces and Helly graphs: most results are presented simultaneously. Moreover, some proofs about injective metric spaces use Helly graphs, and conversely some results about Helly graphs use injective metric spaces. We therefore believe it is worthwile studying these spaces together.

\mk

{\bf Outline of the notes:}  In Section~\ref{sec:definitions}, we start by the various definitions of injective metric spaces and Helly graphs. In Section~\ref{sec:hulls}, we present the crucial tool of injective hulls and Helly hulls. In Section~\ref{sec:local_characterizations}, we present various local-to-global characterizations. In Section~\ref{sec:bicombings_normalforms}, we discuss the possibility to choose nice geodesic paths between any pair of points. In Section~\ref{sec:subdivisions_helly}, we discuss the notion of Helly subdivision of a Helly graph. In Section~\ref{sec:circumcenters_barycenters}, we discuss various notions of centers associated to bounded subsets. In Section~\ref{sec:classification_automorphisms}, we present a classification of automorphisms. In Section~\ref{sec:constructions}, we present various interesting constructions of injective metric spaces and Helly graphs. In Section~\ref{sec:lattices_injective_helly}, we present a very general construction of such spaces from a lattice with a cofinal $1$-parameter action. Towards the end, we gather many known examples and properties of injective and Helly groups (Section~\ref{sec:summary_examples_properties}), and analogies with $L^1$, $L^2$ and $L^\infty$ spaces are sketched (Section~\ref{sec:L1L2Linfini}). Open questions are listed in Section~\ref{sec:questions}, and a few elementary exercises are proposed (gathered in Section~\ref{sec:exercises}).
\mk

{\bf Acknowledgments:} I would like to thank Shaked Bader, Uri Bader, Giuliano Basso, Martin Blufstein, Anthony Genevois, Nima Hoda, Harry Petyt, Urs Lang, Sam Shepherd, Mireille Soergel and Abdul Zalloum for interesting discussions which helped me write these notes. More generally, I would like to thank all the participants of the CRM semester for their enthusiasm, questions and remarks. I would especially like to thank Mark Hagen and Dani Wise for the organization of the CRM semester in Montr\'{e}al and the opportunity to give a minicourse.

\section{Equivalent definitions} \label{sec:definitions}

We will now present the basic and equivalent definitions of injective metric spaces first, and then of Helly graphs.

\subsection{Injective metric spaces}

Our main reference for most definitions and properties of injective metric spaces is~\cite{lang}. In these notes, in a metric space, a ball will always denote a closed ball.

\bdf[Hyperconvex metric space]
A metric space $X$ is called \emph{hyperconvex} if, for any family of points $(x_i)_{i \in I}$ in X, and for any family of radii $(r_i)_{i \in I}$ in $\R_+$ such that for any $i \neq j$ in $I$, we have $r_i+r_j \geq d(x_i,x_j)$, then the intersection $\bigcap_{i \in I} B(x_i,r_i)$ is non-empty.

A metric space $X$ is called \emph{$n$-hyperconvex} (resp. \emph{countably hyperconvex}) if the property holds for any family of at most $n$ balls (resp. for any countable family of balls). 
\edf

\brk
Note that a geodesic metric space $X$ is $3$-hyperconvex if and only if is it \emph{modular}: any three points $x_1,x_2,x_3 \in X$ have a \emph{median} $m \in X$, i.e. for each $i \neq j$, we have $d(x_i,m)+d(m,x_j)=d(x_i,x_j)$. If such a median is unique, the space $X$ is called \emph{metric median}, see for instance~\cite{bowditch_median_metric}.
\erk

\bdf[Injective metric space]
A metric space $X$ is called \emph{injective} if, for any metric space $B$, for any subset $A \subset B$, any $1$-Lipschitz map $f:A \ra X$ has a $1$-Lipschitz extension $\ov{f}:B \ra X$. In other words, $X$ is an injective object in the category of metric spaces, with morphisms corresponding to $1$-Lipschitz maps.
\edf

\bdf[Absolute $1$-Lipschitz retract metric space]
A metric space $X$ is called an \emph{absolute $1$-Lipschitz retract} if, for any isometric embedding $\iota : X \ra Y$, there is a $1$-Lipschitz retraction from $Y$ onto $\iota(X)$.
\edf

\bthm \label{thm:equivalence_definitions_injective}
For any metric space $X$, the following are equivalent:
\bit
\item $X$ is hyperconvex.
\item $X$ is injective.
\item $X$ is an absolute $1$-Lipschitz retract.
\eit
\ethm

\bdf[Injective group] A group is called \emph{(metrically) injective} if it acts properly and cocompactly by isometries on an injective metric space.

A group is called \emph{coarsely injective} if it acts properly and coboundedly by isometries on an injective metric space.
\edf

Before giving a proof of Theorem~\ref{thm:equivalence_definitions_injective}, we will first prove simple results on injective metric spaces.

\blem \label{lem:R_injective}
The space $\R$ is injective.
\elem

\bp
Let $B$ denote a metric space, $A \subset B$ a non-empty subset, and $f:A \ra \R$ a $1$-Lipschitz map. For each $b \in B$, let us define
$$\ov{f}(b) = \sup_{a \in A} (f(a)-d(a,b)).$$
If we fix $a_0 \in A$, for any $a \in A$, we have $f(a)-d(a,b) \leq f(a_0)+d(a_0,a)-d(a,b) \leq f(a_0) + d(a_0,b)$. In particular, the supremum is well-defined.
As a supremum of $1$-Lispchitz maps, it is clear that $\ov{f}$ is $1$-Lipschitz.
Moreover, for any $a,b \in A$, one has $f(a)-d(a,b) \leq f(b)$, hence $\ov{f}(b)=f(b)$. So $\ov{f}$ is the (minimal) $1$-Lipschitz extension of $f$ to $B$. So $\R$ is injective.
\ep

\blem \label{lem:linfinity_product_injective}
Let $(X_i,o_i,d_i)_{i \in I}$ denote a family of pointed injective metric spaces. Their $\ell^\infty$ product is the pointed metric space
$$X=\{(x_i)_{i \in I} \st \forall i \in I, x_i \in X_i \mbox{ and } (d(o_i,x_i))_{i \in I} \mbox{ is bounded.},\}$$
with basepoint $o=(o_i)_{i \in I}$ and metric $d$ defined by:
$$\forall x,y \in X, d(x,y)=\sup_{i \in I} d_i(x_i,y_i).$$
Then $(X,d)$ is injective.
\elem

\bp
Let $(B,d_B)$ denote a metric space, $A \subset B$ a non-empty subset, and $f:A \ra X$ a $1$-Lipschitz map. For each $i \in I$, the $i^\text{th}$ coordinate map $f_i:A \ra X_i$ is $1$-Lipschitz, hence there exists a $1$-Lipschitz extension $\ov{f_i}:B \ra X_i$. For any $b \in B$, $a \in A$, and $i \in I$, we have $d_i(\ov{f_i}(b),o_i) \leq d_i(\ov{f_i}(b),\ov{f_i}(a)) + d_i(f_i(a),o_i) \leq d_B(a,b)+d_i(f_i(a),o_i)$. In particular, the sequence $(d_i(\ov{f_i}(b),o_i))_{i \in I}$ is bounded. The map $\ov{f}=(\ov{f_i})_{i \in I}:B \ra X$ is then a $1$-Lipschitz extension of $f$. So $X$ is injective.
\ep

\brks\

\bit
\item Note that, if $I$ is finite or if each the diameters of the spaces $X_i$ are uniformly bounded, then the choice of basepoints is irrelevant.
\item In particular, if $\Omega$ is any measured space, the metric space $L^\infty(\Omega,\R)$ is injective.
\item In particular, the normed vector space $(\R^n,d_\infty)$, for $n \in \N$, is injective. According to~\cite{nachbin}, these are the only finite-dimensional injective normed vector spaces, up to isometry.
\item As a particular case, since the group $\Z^n$ acts properly an cocompactly by isometries on the injective metric space $(\R^n,d_\infty)$, we deduce that the group $\Z^n$ is injective.
\item More generally, a direct product of injective groups is injective.
\item Note that being injective is not invariant under quasi-isometry, and not even under taking finite index supergroups. For instance, the affine Coxeter group $W \simeq \Z^2 \rtimes \frak{S}_3$ is not coarsely injective (see~\cite{hoda:crystallographic}), but its finite index subgroup $\Z^2$ is injective (and even Helly).
\eit
\erks

\blem \label{lem:retract_is_injective}
A $1$-Lipschitz retract of an injective metric space is injective.
\elem

\bp
Let us consider an injective metric space $Y$, with a subset $X \subset Y$, and a $1$-Lipschitz retract $r:Y \ra X$. Let $B$ denote a metric space, $A \subset B$ a non-empty subset, and $f:A \ra X$ a $1$-Lipschitz map. Since $X \subset Y$ and $Y$ is injective, $f$ has a $1$-Lipschitz extension $\ov{f}:B \ra Y$. Now $r \circ \ov{f} : B \ra X$ is a $1$-Lipschitz extension, and $X$ is injective.
\ep

\bp[Proof of Theorem~\ref{thm:equivalence_definitions_injective}]
Assume that $X$ is injective, we will prove that $X$ is hyperconvex. Consider a family of points $(x_i)_{i \in I}$ in X, and for any family of radii $(r_i)_{i \in I}$ in $\R_+$ such that for any $i \neq j$ in $I$, we have $r_i+r_j \geq d(x_i,x_j)$. Note that we may assume that the points $(x_i)_{i \in I}$ are pairwise distinct. Consider the metric space $A$ with underlying set $I$, such that for all $i,j \in A$, we set $d_A(i,j)=d(x_i,x_j)$. Consider the metric space $B=A \cup \{y\}$, containing an isometric copy of $A$, where for all $i \in A$, we define $d_B(i,y) = r_i$.  The natural map $f:i \in A \mapsto x_i \in X$ is an isometric embedding, so since $X$ is injective it extends to a $1$-Lipschitz $\ov{f} : B \ra X$. The point $z=\ov{f}(y) \in X$ is such that, for each $i \in I$, we have $d(z,x_i) \leq d_B(y,i) = r_i$. In particular, the point $z$ lies in the intersection of all balls $B(x_i,r_i)$. So $X$ is hyperconvex.

\mk

Assume that $X$ is hyperconvex, we will prove that $X$ is an absolute $1$-Lipschitz retract. Consider an isometric embedding $\iota:X \ra Y$ into a metric space $Y$: to simplify notations, assume that $X \subset Y$. Consider the set ${\cal Z}$ of subsets $Z \subset Y$ such that $X \subset Z$ and $X$ is a $1$-Lipschitz retract of $Z$. This is an inductive family for the reverse inclusion order, so according to Zorn's Lemma, we may consider a maximal element $Z \in {\cal Z}$, and let $r:Z \ra X$ denote a $1$-Lipschitz retraction. Assume by contradiction that $Z \subsetneq Y$, and let $z \in Y \bs Z$. For each $x,y \in Z$, we know by the triangle inequality that $d(x,z)+d(y,z) \geq d(x,y) \geq d(r(x),r(y))$. Since $X$ is hyperconvex, there exists a point $\ov{z} \in X$ such that, for any $x \in X$, we have $\ov{z} \in B(r(x),d(x,z))$. In particular, if extend $r:Z \ra X$ to $\ov{r}:Z \cup \{z\} \ra X$ by $\ov{r}(z)=\ov{z}$, the map $\ov{r}$ is a $1$-Lipschitz retract, which contradicts the maximality of $Z$. Hence $Z=Y$, and $X$ is an absolute $1$-Lipschitz retract.

\mk

Assume that $X$ is an absolute $1$-Lipschitz retract. Let $Y$ denote the space of bounded functions $X \ra \R$, endowed with the $d_\infty$ metric. Fix $x_0 \in X$, and let us consider the following isometric embedding
\beq \iota:X & \ra & Y \\
x & \mapsto & d(x,\cdot)-d(x_0,\cdot).\eeq
According to Lemma~\ref{lem:linfinity_product_injective}, the space $Y$ is injective. Since $X$ is an absolute $1$-Lipschitz retract, we deduce that there exists a $1$-Lipschitz retract from $Y$ to $\iota(X)$. According to Lemma~\ref{lem:retract_is_injective}, we deduce that $X$ is injective.
\ep

\blem \label{lem:injective_implies_geodesic}
Any injective metric space is geodesic.
\elem

\bp
Consider an injective metric space $X$, and two distinct points $x,y \in X$ at distance $r=d(x,y)$. Consider the interval $B=[0,r]$, and its boundary $A=\{0,r\}$. Consider the isometric embedding $f:A \ra X$ defined by $f(0)=x$ and $f(r)=y$. Since $X$ is injective, $f$ extends to a $1$-Lipschitz map $\ov{f}:[0,r] \ra X$. Since $d(f(0),f(r))=d(x,y)=r$, we deduce that $\ov{f}$ is an isometric embedding, i.e. it is a geodesic from $x$ to $y$. So $X$ is geodesic.
\ep

\blem \label{lem:injective_implies_complete}
Any injective metric space is complete.
\elem

\bp
Consider an injective metric space $X$. Let $\iota:X \ra \ov{X}$ denote the embedding into the metric completion of $X$. Since $X$ is an absolute $1$-Lipschitz retract according to Theorem~\ref{thm:equivalence_definitions_injective}, we deduce that there exists a $1$-Lipschitz retract from $\ov{X}$ to $\iota(X)$. As $\iota(X)$ is dense in $\ov{X}$, we conclude that $\iota(X)=\ov{X}$, i.e. $X$ is complete.
\ep

In the case of geodesic metric spaces, we deduce the simplest way to define injective metric spaces. Let us recall first a definition.

\bdf[Helly property]
A family ${\cal F}$ of subsets of a set $X$ is said to have the \emph{Helly property} if any family of pairwise intersecting elements of ${\cal F}$ has a non-empty global intersection.
\edf

\bthm \label{thm:injective_equivalent_pairwise_intersecting_balls}
A geodesic metric space $X$ is injective if and only if (closed) balls satisfy the Helly property.
\ethm

\bp
We will prove that any such space $X$ is hyperconvex. Indeed, if $x_i,x_j \in X$ and $r_i,r_i \geq 0$ are such that $d(x_i,x_j) \leq r_i+r_j$, then since $X$ is geodesic the balls $B(x_i,r_i)$ and $B(x_j,r_j)$ intersect.
\ep

\brk
One could also wonder about injective objects in the category of $\Lambda$-valued metric spaces, for some abelian group $\Lambda$. For instance, one could consider Weyl-chamber-valued metric spaces as for buildings and symmetric spaces, see~\cite{parreau_compactifications}.
\erk

Let us remark that there are plenty of examples of injective spaces other than mere $\ell^\infty$ normed vector spaces.

\bpro
Any complete $\R$-tree is injective.
\epro

\bp
Let $X$ denote a complete $\R$-tree. It is geodesic, so according to Theorem~\ref{thm:injective_equivalent_pairwise_intersecting_balls}, it is sufficient to prove that balls in $X$ satisfy the Helly property. 

\mk

We will first prove that $X$ is $3$-hyperconvex: let us consider three pairwise intersecting balls $(B(x_i,r_i))_{1 \leq i \leq 3}$ in $X$. Let $m \in X$ denote the median of $x_1,x_2,x_3$, i.e. the unique point such that $[x_1,x_2] \cap [x_2,x_3] \cap [x_3,x_1] = \{m\}$. If $r_i \geq d(m,x_i)$ for all $i \in I$, then $m \in \cap_{i \in I} B(x_i,r_i)$. Assume that there exists $i \in I$, say $i=1$, such that $r_1 < d(m,x_1)$. Then the point $y$ on $[x_1,m]$ at distance $r_1$ from $x_1$ is such that $y \in \cap_{i \in I} B(x_i,r_i)$. Hence $X$ is $3$-hyperconvex.

\mk

Now consider an arbitrary family $(B(x_i,r_i))_{i \in I}$ of (at least $3$) pairwise intersecting balls of $X$. Fix $i_0 \in I$, then we may assume up to decreasing $r_{i_0}$ that we have
$$r_{i_0} = \sup_{i \in I} d(x_{i_0},x_i)-r_i.$$
In particular, for each $n \geq 1$, there exists $i_n \in I$ such that $d(x_{i_0},x_{i_n}) \geq r_{i_0}+r_{i_n}-\f{1}{2^n}$.

For each $n \geq 1$, the intersection $A_n=B(x_{i_0},r_{i_0}) \cap B(x_{i_n},r_{i_n})$ has diameter at most $\f{2}{2^n}$. For each $n \geq 1$, let us consider $y_n \in A_n \cap A_{n+1}$, which is non-empty according to $3$-hyperconvexity. Since $y_n,y_{n+1} \in A_{n+1}$, we deduce that $d(y_n,y_{n+1}) \leq \f{1}{2^n}$. So the sequence $(y_n)_{n \in \N}$ is Cauchy: let us denote its limit by $y \in X$.

We will now prove that, for each $i \in I$, we have $d(x_i,y) \leq r_i$. By contradiction, assume that there exists $i \in I$ such that $d(x_i,y) > r_i$. Let $n \geq 1$ such that $d(x_i,y) > r_i+\f{2}{2^n}$. Since $X$ is $3$-hyperconvex, there exists $z \in A_n \cap B(x_i,r_i)$. We deduce that $d(y,z)>\f{2}{2^n}$, which contradicts $y,z \in A_n$ and $\diam(A_n) \leq \f{2}{2^n}$. Hence $y \in \bigcap_{i \in I} B(x_i,r_i)$: $X$ is hyperconvex, so it is injective.
\ep

We deduce immediately the following.

\bcor
Any product of complete $\R$-trees, endowed with the sup metric, is injective.
\ecor

\bcor
Any finite rank free group is injective.
\ecor

\subsection{Helly graphs}

Helly graphs are the perfect discrete counterpart of injective metric spaces. Rather than being a parallel theory, there is a rich interplay between injective metric spaces and Helly graphs, as we will see. Good references for Helly graphs are~\cite{chalopin_chepoi_hirai_osajda} and~\cite{helly_groups}.

\bdf[Discretely geodesic]
A metric space $X$ with integer-valued metric is called \emph{discretely geodesic} if, for any $x,y \in X$ with $d(x,y)=n$, there exist $x_0=x,x_1,\dots,x_n=y$ in $X$ such that, for each $0 \leq i \leq n-1$, we have $d(x_i,x_{i+1})=1$.
\edf

\brk
A metric space $X$ is discretely geodesic if and only if $X$ is the vertex set of a connected graph, endowed with the combinatorial distance. As a slight abuse of notation, we will often identify a connected graph with its vertex set, endowed with the combinatorial distance.
\erk

\bdf[Helly graph]
A connected graph $X$ is called \emph{Helly} if any family of pairwise intersecting combinatorial balls has a non-empty intersection.
\edf

\bdf[Helly group] A group is called \emph{Helly} if it acts properly cocompactly by automorphisms on a Helly graph.

A group is called \emph{coarsely Helly} if it acts properly coboundedly by automorphisms on a Helly graph.
\edf

We will see in Corollary~\ref{cor:helly_implies_injective} that every Helly group is an injective group, and every coarsely Helly group is a coarsely injective group.

\mk

Instead of working in the category of metric spaces, we can restrict ourselves to the subcategory of metric spaces with integer-valued metric. We will use the adjective integral to refer to the corresponding properties of being injective, absolute $1$-Lipschitz
retract or hyperconvex. As for injective metric spaces, we have the following equivalent characterizations of Helly graphs.

\bthm \label{thm:equivalent_characterizations_Helly_graphs}
Consider a metric space $X$ with integer-valued metric. Then the following are equivalent:
\bit
\item $X$ is integrally injective.
\item $X$ is an absolute $1$-Lipschitz integral retract.
\item $X$ is integrally hyperconvex.
\item $X$ is the vertex set of a Helly graph, with the combinatorial metric.
\eit
\ethm

\bpro
Any simplicial tree is a Helly graph.
\epro

\bp
Let us consider a family of pairwise intersection balls $(B(x_i,r_i))_{i \in I}$ in a tree $X$. Up to reducing radii, we may assume that there exist $i,j \in I$ such that $d(x_i,x_j)=r_i+r_j$. Hence $B(x_i,r_i) \cap B(x_j,r_j)=\{z\}$. Now, for each $k \in I \bs \{i,j\}$, since $B(x_k,r_k)$ is connected and intersects $B(x_i,r_i)$ and $B(x_j,r_j)$, we deduce that $z \in B(x_k,r_k)$ as $X$ is a tree. Hence $z \in \bigcap_{i \in I} B(x_i,r_i)$: the tree $X$ is Helly.
\ep

\bexe (Exercise)
Fix $n \geq 1$. Consider the graph $\Gamma$ with vertex set $\Z^n$, with an edge between $v$ and $w$ if, for all $1 \leq i \leq n$, we have $|v_i-w_i| \leq 1$. Then $\Gamma$ is a Helly graph.
\eexe

More generally, we may consider products of Helly graphs: the proof is similar to that of Lemma~\ref{lem:linfinity_product_injective}.

\blem \label{lem:linfinity_product_Helly}
Let $(X_i,o_i)_{i \in I}$ denote a family of pointed Helly graphs, where $o_i$ is a vertex of $X_i$ for each $i \in I$. Their $\ell^\infty$ product is the graph $X$ with vertex set
$$V(X)=\{(x_i)_{i \in I} \st \forall i \in I, x_i \in X_i \mbox{ and } (d(o_i,x_i))_{i \in I} \mbox{ is bounded.},\}$$
with an edge between $(x_i)_{i \in I}$ and $(y_i)_{i \in I}$ if, for each $i \in I$, either $x_i=y_i$ or $x_i$ is adjacent to $y_i$.
Then $X$ is a Helly graph.
\elem

\bexes We can therefore find elementary examples of Helly groups:
\bit
\item For any $n \geq 0$, the free abelian group $\Z^n$ is Helly.
\item For any $n \geq 0$, the free group $\F_n$ is Helly.
\eit
\eexes

\section{Hulls} \label{sec:hulls}

\subsection{Injective hulls}

A key notion in the theory of injective metric spaces is that of an injective hull. This has been discovered by Isbell in 1963 (see~\cite{isbell,holsztynski_injective_hull}) and rediscovered later (see~\cite{dress,dress_moulton_terhalle_Ttheory,chrobak_larmore_injective_hull_algorithm,holsztynski_injective_hull}). This construction has several names: "injective hull", "injective envelope", "hyperconvex hull" or "tight span".

\bthm[\cite{isbell}] \label{thm:existence_injective_hull}
For any metric space $X$, there is an injective metric space $E(X)$ called the \emph{injective hull} of $X$, and an isometric embedding $e:X \ra E(X)$, which is minimal in the following sense. For any isometric embedding $f:X \ra Y$ of $X$ into an injective metric space $Y$, the embedding $f$ factors by $e$, i.e. there exists an isometric embedding $E(f):E(X) \ra Y$ such that $f=E(f) \circ e$.

Moreover, $E(X)$ is essentially unique, in the following sense. If $e:X \ra E$ and $e':X \ra E'$ are two injective hulls, there exists a unique isometry $\iota:E \ra E'$ such that $e'=\iota \circ e$.
\ethm

Here is an explicit description of the injective hull of any metric space (see~~\cite{dress,dress_moulton_terhalle_Ttheory,chrobak_larmore_injective_hull_algorithm,holsztynski_injective_hull,lang}).

Let $(X,d)$ denote an arbitrary metric space, and let $\Lip_1(X,\R)$ denote the vector space of $1$-Lipschitz maps from $X$ to $\R$. Let us define
$$\Delta(X)=\{f \in \Lip_1(X,\R) \st \forall x,y \in X, f(x) + f(y) \geq d(x,y)\}.$$
Let us consider the sup metric $d_\infty$ on $\Delta(X)$. We then have the canonical Kuratowski isometric embedding
\beq e:(X,d) & \mapsto & (\Delta(X),d_\infty) \\
x & \mapsto & (y \mapsto d(x,y)).\eeq

\bthm \label{thm:explicit_description_injective_hull}
The subspace
\beq E(X)&=&\{f \in \Delta(X) \mbox{ minimal}\} \\
&=&\{f \in \Delta(X) \st \forall g \in \Delta(X), g \leq f \Rightarrow g=f\} \\
&=&\{f \in \R^X \st \forall x \in X, f(x) = \sup_{y \in X} (d(x,y)-f(y))\}\eeq
with the sup metric $d_\infty$ and the isometric embedding $e:(X,d) \ra (E(X),d_\infty)$, is the injective hull of $X$.
\ethm

We will follow Lang's proof to this result, from which we state here the main steps.

\blem \label{lem:delta_injective}
The space $(\Delta(X),d_\infty)$ is injective.
\elem

\bp
We give here a direct proof of this result.

Let us consider a family $(f_i)_{i \in I}$ in $\Delta(X)$, and a family $(r_i)_{i \in I}$ in $\R_+$, such that $\forall i,j \in I, d_\infty(f_i,f_j) \leq r_i +r_j$. Let us define
\beq g : X & \ra & \R_+\\
x & \mapsto & \inf_{i \in I} f_i(x)+r_i.\eeq

We will prove that $g \in \Delta(X)$.

Let us first prove that $g$ is $1$-Lipschitz. Fix $x,y \in X$, $\eps>0$, and let $i \in I$ such that $g(x) \geq f_i(x)+r_i-\eps$. Then we have
$$g(x)-g(y) \geq f_i(x)+r_i-\eps-f_i(y)-r_i \geq -d(x,y)-\eps.$$
Since this holds for any $\eps>0$, we deduce that $g(x)-g(y) \geq -d(x,y)$, hence $g$ is $1$-Lipschitz.

For each $x,y \in X$, let $\eps>0$, and let $i,j \in I$ such that $g(x) \geq f_i(x)+r_i-\eps$ and $g(y) \geq f_j(y)+r_j-\eps$. We deduce that
$$g(x)+g(y) \geq f_i(x)+r_i-\eps+f_j(y)+r_j-\eps \geq f_j(x)+f_j(y) -2\eps \geq d(x,y)-2\eps.$$
Since this holds for any $\eps>0$, we deduce that $g(x)+g(y) \geq d(x,y)$. We conclude that $g \in \Delta(X)$.

\mk

Let us finally prove that, for each $i \in I$, we have $d_\infty(g,f_i) \leq r_i$. Fix $x \in X$, we first have $g(x) \leq f_i(x)+r_i$. Moreover, for any $j \in I$ we have $f_j(x) \geq f_i(x)-r_i-r_j$, hence $g(x) \geq f_i(x)-r_i$. So we deduce that $d_\infty(g,f_i) \leq r_i$.

So we have proved that $\Delta(X)$ is hyperconvex. According to Theorem~\ref{thm:equivalence_definitions_injective}, this implies that $\Delta(X)$ is injective.
\ep

\blem \label{lem:delta_retract}
There is a $1$-Lipschitz retraction $p:\Delta(X) \ra E(X)$, which is equivariant with respect to the isometry group of $X$.
\elem

\bp
We follow here Lang's proof of~\cite[Proposition~3.1]{lang}.

For each $f \in \Delta(X)$, let us define
\beq f^\star : X & \ra & \R \\
x & \mapsto & \sup_{z \in X} d(x,z)-f(z),\eeq
and let $q(f)=\f{1}{2}(f+f^\star)$: for each $x,y \in X$, we have $f(x)+f^\star(y)\geq d(x,y)$ and $f(y)+f^\star(x)\geq d(x,y)$, hence $q(f)(x)+q(f)(y) \geq d(x,y)$.

Moreover, we will see that $f^\star$ is $1$-Lipschitz: fix $x,y \in X$ and $\eps>0$, and assume that $z \in X$ is such that $f^\star(x) \leq d(x,z)-f(z)+\eps$: since $f^\star(y) \geq d(y,z)-f(z)$, we deduce that $f^\star(x)-f^\star(y) \leq d(x,z)-d(y,z)+\eps \leq d(x,y)+\eps$. Since this holds for $\eps>0$, by symmetry, we deduce that $f^\star$ is $1$-Lipschitz. Hence $q(f)$ is also $1$-Lipschitz, and $q(f) \in \Delta(X)$.

\mk

Furthermore, remark that for any $f \in \Delta(X)$, for all $x,z \in X$ we have $f(x)+f(z) \geq d(x,z)$, hence $f^\star(x) \leq f(x)$. In particular, $q(f) \leq f$. We may thus define $p(f):X \ra \R$ as the pointwise limit of the non-increasing sequence of non-negative functions $(q^n(f))_{n \in \N}$. It is clear that $p(f) \in \Delta(X)$, let us prove that $p$ is a $1$-Lipschitz retraction onto $E(X)$.

\mk

Note that, for any $f,g \in \Delta(X)$, we have $d_\infty(f^\star,g^\star) \leq d_\infty(f,g)$, so we have $d_\infty(q(f),q(g)) \leq d_\infty(f,g)$: we deduce that $p$ is $1$-Lipschitz.

\mk

Now remark that $E(X)$ is precisely the fixed point set of $p$. Furthermore, fix $f \in \Delta(X)$, we will show that $p(f) \in E(X)$. Remark that, for each $n \geq 1$, we have $p(f) \leq q^n(f)$, hence $p(f)^\star \geq q^n(f)^\star$, so
$$0 \leq p(f)-p(f)^\star \leq q^n(f)-q^n(f)^\star \leq 2(q^n(f)-q^{n+1}(f)).$$
Since the sequence of functions $(2(q^n(f)-q^{n+1}(f))_{n \in \N}$ converges to $0$ as $n \ra \infty$, we conclude that $p(f)=p(f)^\star$, and so $p(f) \in E(X)$. Hence $p$ is a $1$-Lipschitz retract on $E(X)$.

It is clear by definition that $p$ is equivariant with respect to the isometry group of $X$.
\ep

We can now finish the proof of Theorem~\ref{thm:explicit_description_injective_hull}.

\bp[of Theorem~\ref{thm:explicit_description_injective_hull}]
According to Lemma~\ref{lem:delta_injective}, we know that $\Delta(X)$ is injective. According to Lemma~\ref{lem:delta_retract}, we know that $E(X)$ is a $1$-Lipschitz retract of $\Delta(X)$. According to Lemma~\ref{lem:retract_is_injective}, we conclude that $E(X)$ is injective.

\mk

Let us now prove that $E(X)$ is an injective hull of $X$: let us assume that $f:X \ra Y$ is an isometric embedding into an injective metric space $Y$. Since $Y$ is injective, there exists a $1$-Lipschitz map $\phi:E(X) \ra Y$ such that $\phi \circ e=f$. Since $E(X)$ is injective, there exists a $1$-Lipschitz map $\psi:Y \ra E(X)$ such that $\psi \circ f=e$. Hence $\theta=\psi \circ \phi:E(X) \ra E(X)$ is a $1$-Lipschitz map which restricts to the identity on $e(X)$: we will prove that $\theta=\id$. 

By contradiction, assume that there exists $g \in E(X)$ such that $\theta(g) \neq g$. For any $x \in X$, we have
\beq \theta(g)(x) &=& d(\theta(g),e(x)) \\
&=& d(\theta(g),\theta(e(x))) \\
& \leq & d(g,e(x)) = g(x).\eeq
So we have $\theta(g) \leq g$: since $g$ is minimal in $\Delta(X)$, we conclude that $\theta(g)=g$, so $\theta=\id$.
\ep

One immediate consequence of the proof is the following.

\bcor \label{cor:injective_contractible}
Any injective metric space is contractible.
\ecor

\bp
Let $X$ denote an injective metric space. The subspace $\Delta(X)$ of $\R^X$ is affinely convex, hence it is contractible. Since $p:\Delta(X) \ra E(X) \simeq X$ is a retraction, we conclude that $E(X)$ is contractible.
\ep

\bexes\

\begin{figure}[H]
\begin{center}
\begin{tikzpicture}
\def \p {0.05}
\def \op {0.5}
\def \gris {black!10}
\draw[fill] (1,0) circle (\p) node(A) {};
\draw[fill] (0,1) circle (\p) node(B) {};
\draw[fill] (4,3) circle (\p) node(D) {};
\draw[fill] (3,4) circle (\p) node(C) {};
\draw[fill] (1,-1) circle (\p) node(a) {};
\draw[fill] (-2,1) circle (\p) node(b) {};
\draw[fill] (5.5,3) circle (\p) node(d) {};
\draw[fill] (3,5) circle (\p) node(c) {};

\draw[black,fill opacity=\op,fill=\gris] (A.center) -- (B.center) -- (C.center) -- (D.center) -- cycle;
\draw (A.center) -- (a.center);
\draw (B.center) -- (b.center);
\draw (C.center) -- (c.center);
\draw (D.center) -- (d.center);
\node at (2,2) {\bfseries $\ell^\infty$};

\end{tikzpicture}
\end{center}
\caption{The injective hull of a $4$-point metric space}
\label{fig:hull_4point}
\end{figure}
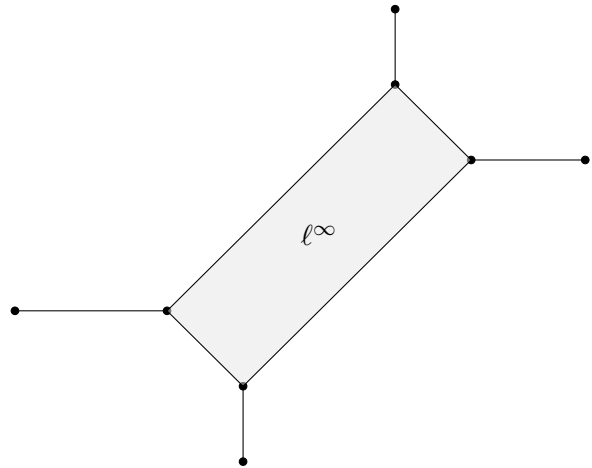

\bit
\item If $X$ is already injective, then $E(X)=X$.
\item If $X$ is a metric space consisting of $2$ points $\{a,b\}$, then $E(X)$ is the segment of length $d(a,b)$.
\item If $X$ is a metric space consisting of $3$ points, then $E(X)$ is a tripod (Exercise).
\item If $X$ is a metric space consisting of $4$ points, then $E(X)$ is a "slanted rectangle" in $(\R^2,\ell^\infty)$ with antennas attached to the four corners (see Figure~\ref{fig:hull_4point}).
\item If $X$ is a metric space consisting of $5$ points, there are three generic combinatorial types, which do not look right-angled anymore, see~\cite{dress}.
\eit
\eexes

\subsection{Helly hulls}

We also have the existence of a minimal Helly graph containing isometrically any given connected graph.

\bthm \label{thm:exists_helly_hull}
For any simplicial connected graph $X$, there is a Helly graph $H(X)$ called the \emph{Helly hull} of $X$, and an isometric simplicial embedding $e:X \ra H(X)$, which is minimal in the following sense. For any isometric simplicial embedding $f:X \ra Y$ of $X$ into a Helly graph $Y$, the embedding $f$ factors by $e$, i.e. there exists an isometric simplicial embedding $H(f):H(X) \ra Y$ such that $f=H(f) \circ e$.

Moreover, $H(X)$ is essentially unique, in the following sense. If $e:X \ra Y$ and $e':X \ra Y'$ are two Helly hulls, there exists a unique simplicial isomorphism $\iota:Y \ra Y'$ such that $e'=\iota \circ e$.
\ethm

\brk
In particular, if $G$ is a group of automorphisms of a connected graph $X$ (or a group of isometries of a metric space $Y$), then the action of $G$ extends uniquely to an action by automorphisms of the Helly hull $H(X)$ (or by isometries on the injective hull $E(Y)$). In particular, any finitely generated group has a metrically proper action on a Helly graph, the Helly hull of any Cayley graph.

This explains why we often place restrictions on actions on Helly graphs and injective metric spaces, like coboundedness assumptions or dimension bounds.
\erk

As in the case of the injective hull of a metric space, there is an explicit model for the Helly hull. Let us consider the vertex set $X$ of a connected graph. One can essentially think of the Helly hull of $X$ as being the $1$-skeleton of the injective hull of $X$, in the following more precise way. Let us define
$$\Delta(X)=\{f \in \N^X \mbox{ $1$-Lipschitz} \st \forall x,y \in X, f(x)+f(y) \geq d(x,y)\}.$$

\bthm \cite[Theorem~4.4]{helly_groups} \label{thm:helly_hull_integer_points_inject_hull}
Let $X$ denote the vertex set of a connected graph, let $H(X)$ denote the Helly hull of $X$, and let $E(X)$ denote the injective hull of $X$ endowed with the combinatorial distance. Then the vertex set of $H(X)$ is
\beq H(X)^{(0)} &=& E(X) \cap \N^X\\
&=& \{f \in \Delta(X) \mbox{ minimal}\} \\
&=&\{f \in \Delta(X) \st \forall g \in \Delta(X), g \leq f \Rightarrow g=f\} \\
&=&\{f \in \N^X \st \forall x \in X, f(x) = \max_{y \in X} (d(x,y)-f(y))\},\eeq
with an edge between $f,g \in E(X) \cap \N^X$ if and only if $d_\infty(f,g) =1$.
\ethm

\bexe
This characterization is sufficient to compute the Helly hull of some small examples. For instance, let us compute the Helly hull of a $4$-cycle $X$ with vertices $(a_i)_{i \in \Z/4\Z}$ in this order. Let $f:X \ra \N$ denote a minimal element of $\Delta(X)$.

Assume first that there exists $i \in \Z/4\Z$ such that $f(a_i)=0$. Then, since $f$ belongs to $\Delta(X)$, we see that $f \geq e(a_i)$, hence by minimality we have $f=e(a_i)$.

Assume now that there exists $i \in \Z/4\Z$ such that $f(a_i)=2$. Then it is not hard to see that $f \geq e(a_{i+2})$, hence by minimality we have $f=e(a_{i+2})$.

We are left with the case $\forall i \in \Z/4\Z, f(a_i)=1$. This defines a minimal element $f$ of $\Delta(X)$. Hence the Helly hull $H(X)$ is the cone over the $4$-cycle, with cone vertex $f$.
\eexe

\bexe (Exercise)
\bit
\item What is the Helly hull of a $5$-cycle?
\item What is the Helly hull of a $6$-cycle?
\eit
\eexe

We will study in Section~\ref{sec:subdivisions_helly} more precisely a cell structure on the injective hull of a graph.

\mk

One of the most powerful tools to control the geometry of hulls of graphs follows from the work of Lang.

\bdf[Stable intervals]
A graph $X$ is said to have \emph{stable intervals} if there exists $\beta \geq 1$ such that, for any $x,y,z \in X$ such that $d(y,z)=1$, for any geodesic $[x,y]$, there exists a geodesic $[x,z]$ at Hausdorff distance at most $\beta$ from $[x,y]$.
\edf

\bexe (Exercise)\
\ben
\item Any Gromov-hyperbolic graph has stable intervals.
\item Any Helly graph has stable intervals, with constant $\beta=1$.
\item Any median graph (i.e. $1$-skeleton of a CAT(0) cube complex) has stable intervals, with constant $\beta=1$.
\een
\eexe

The following result by Lang is crucial in the study of injective and Helly hulls.

\bthm \cite[Theorem~1.1]{lang} \label{thm:stable_intervals_proper_hull}
Let $X$ denote a locally finite graph with stable intervals. Then $H(X)$ is locally finite, and $E(X)$ is proper, and locally finite-dimensional.
\ethm

\section{Simple characterizations} \label{sec:local_characterizations}

We will now see simple and powerful tools to characterize injective metric spaces. Notably, every notion of (nonpositive) curvature should be local. Hence we expect to have local-to-global results for injective metric spaces and Helly graphs.

\subsection{Local characterizations of Helly graphs}

If $X$ is a graph, its \emph{triangle complex} is the $2$-dimensional simplicial complex whose $2$-simplices are triangles of $X$. Recall that a \emph{clique} of $X$ is a complete subgraph. The following is a deep result by Chalopin, Chepoi, Hirai and Osajda (see~\cite{chalopin_chepoi_hirai_osajda}). The last condition of clique-Helly is usually the simplest to check.

\bthm[Local-to-global for Helly graphs] \label{thm:local_global_Helly}
Let $X$ denote a connected graph. Then the following are equivalent:
\bit
\item $X$ is Helly.
\item {\bf 1-Helly} The triangle complex of $X$ is simply connected, and $1$-balls in $X$ satisfy the Helly property.
\item {\bf Clique-Helly} The triangle complex of $X$ is simply connected, and maximal cliques in $X$ satisfy the Helly property.
\eit
\ethm

\subsection{Cartan-Hadamard theorem for injective metric spaces}

As for every fine notion of nonpositive curvature, there is a local-to-global result (the analogue of the Cartan-Hadamard Theorem for CAT(0) spaces). Say that a metric space is \emph{uniformly locally injective} if there exists $r>0$ such that each ball of radius $r$ is injective.

Note that, according to Corollary~\ref{cor:injective_contractible}, any injective metric space is simply connected.

\bthm[Local-to-global for injective metric spaces, \cite{miesch,haettel_helly_kpi1}] \label{thm:local_global_injective_spaces}
Let $X$ denote a complete, simply connected, uniformly locally injective metric space. Then $X$ is injective.
\ethm

Note that Miesch's proof (see~\cite{miesch}) requires properness, and follows the proof of the Cartan-Hadamard Theorem from~\cite{bridson_haefliger} using bicombings. However, our proof (see the main argument below, and~\cite{haettel_helly_kpi1} for the complete details) does not require properness, but relies on the local-to-global property for Helly graphs (Theorem~\ref{thm:local_global_Helly}). 

\bp
Fix $\eps>0$ small enough such that balls in $X$ of radius at most $2\eps$ are injective. Consider the graph $\Gamma_\eps$ with vertex set $X$, and with an edge between $x,y \in X$ if $d(x,y) \leq \eps$. Since $X$ is path-connected, $\Gamma_\eps$ is a connected graph. Also note that, for any $x \in X$ and $n \in \N$, we have
$$B_{\Gamma_\eps}(x,n) = B_X(x,n\eps).$$

\mk

We will prove that, for each $\eps>0$, the graph $\Gamma_\eps$ is a Helly graph by applying Theorem~\ref{thm:local_global_Helly}: more precisely, we will prove that $\Gamma_\eps$ is $1$-Helly.

\mk

Fix a family of vertices $(x_i)_{i \in I}$ of $\Gamma_\eps$ such that $\forall i,j \in I, d_{\Gamma_\eps}( x_i,x_j) \leq 2$. We want to prove that the balls $(B_{\Gamma_\eps}(x_i,1))_{i \in I}$ intersect in $\Gamma_\eps$.

The family of metric balls $(B_X(x_i,\eps))_{i \in I}$ in $X$ pairwise intersects: since such balls have the Helly property by assumption on $X$, so we deduce that there exists $y \in X$ such that $\forall i \in I, d_X(x_i,y) \leq \eps$. In other words, the vertex $y \in \Gamma_\eps$ lies in the intersection of all combinatorial $1$-balls $(B_{\Gamma_\eps}(x_i,1))_{i \in I}$. We deduce that the graph $\Gamma_\eps$ is $1$-Helly.

\mk

We now prove that the triangle complex of $\Gamma_\eps$ is simply connected. Fix a combinatorial loop $\ell$ in $\Gamma_\eps$. Since $X$ is simply connected, there exists a disk $D$ in $X$ bounding $\ell$. Consider a triangulation $T$ of $D$ such that triangles have diameter for $d_X$ at most $\eps$. Then the vertex set of each triangle of $T$ is a clique in $\Gamma_\eps$, therefore $\ell$ is null-homotopic in the triangle complex of $\Gamma_\eps$. So the triangle complex of $\Gamma_\eps$ is simply connected.

\mk

According to Theorem~\ref{thm:local_global_Helly}, we deduce that $\Gamma_\eps$ is Helly. This implies that the metric space $X$ is $\eps$-coarsely injective (see Section~\ref{sec:coarse_injectivity}), for each $\eps>0$. As $X$ is complete, Proposition~\ref{pro:epsilon_coarse_injective} implies that $X$ is injective.
\ep

\subsection{$4$-hyperconvexity}

For $n \geq 3$, recall that a metric space $(X,d)$ is $n$-hyperconvex if the property hyperconvexity holds for $n$ balls: for any family $(x_i)_{1 \leq i \leq n} \in X^n$, and and family $(r_i)_{1 \leq i \leq n} \in \R_+^n$ such that for all $1 \leq i,j \leq n$ we have $r_i+r_j \geq d(x_i,x_j)$, then $\bigcap_{i=1}^n B(x_i,r_i) \neq \emptyset$.

\mk

Miesch and Pavon proved that for a complete metric space, $4$-hyperconvexity implies finite hyperconvexity (see~\cite{miesch_pavon_4hyperconvex}). With Hoda and Petyt, we extended this result to hyperconvexity with some extra assumptions.

\bthm Miesch-Pavon~\cite{miesch_pavon_4hyperconvex}, Haettel-Hoda-Petyt~\cite[Work in progress]{haettel_hoda_petyt_4hyperconvexity} \label{thm:4hyperconvex_implies_injective}
Let $X$ denote a proper metric space. If $X$ is $4$-hyperconvex, then $X$ is injective.
\ethm

\bp
Since $X$ is proper, it is complete. According to~\cite{miesch_pavon_4hyperconvex}, we deduce that $X$ is $n$-hyperconvex for all $n \geq 3$.

We will now prove that $X$ is countably hyperconvex. Let $(B_k)_{k \in \N}$ denote a countable family of pairwise intersecting balls of $X$. For each $k \in \N$, the intersection $C_k=\bigcap_{h \leq k} B_h$ is non-empty by $(k+1)$-hyperconvexity: the decreasing sequence of compact subsets $(C_k)_{k \in \N}$ has non-empty global intersection, so the balls $(B_k)_{k \in \N}$ have a non-empty global intersection: $X$ is countably hyperconvex.

Since $X$ is proper, it is separable. As remarked by Hoda in~\cite{hoda:crystallographic}, $X$ is Lindel\"of: for any family ${\cal F}$ of closed subsets of $X$ such that any countable subfamily has nonempty intersection, the family ${\cal F}$ itself has nonempty intersection. In particular, $X$ is hyperconvex.
\ep

One notable consequence of Theorem~\ref{thm:4hyperconvex_implies_injective} is the following.

\bcor \label{cor:asymptotic_cone_injective}
Let $X$ denote a $4$-hyperconvex metric space, and assume that an asymptotic cone $Y$ of $X$ is proper. Then $Y$ is injective.
\ecor

\bp
It is not hard to see that any asymptotic cone $Y$ of $X$ is also $4$-hyperconvex. Also, any asymtptotic cone is complete. According to Theorem~\ref{thm:4hyperconvex_implies_injective}, we deduce that $Y$ is injective.
\ep

Note that the asymptotic cone of a finitely generated group $G$ is proper if and only if $G$ is virtually nilpotent (this is essentially due to Gromov, see~\cite[Remark~2.8]{sapir_lecture_asymptotic}).

One consequence is a very simple example of a group with no proper cobounded action on an injective metric spaces. For the family of crystallographic groups, see~\cite{hoda:crystallographic}.

\bcor[Hoda]
Let $W \simeq \Z^2 \rtimes \frak{S}_3$ denote the $(3,3,3)$ triangle Coxeter group, i.e. the group generated by Euclidean reflections with respect to the lines of the standard equilateral tiling of $\R^2$. Then $W$ has no proper cobounded action on an injective metric space.
\ecor

\bp
By contradiction, assume that $W$ acts properly and coboundedly on an injective metric space $X$. Let $Y$ denote an asymptotic cone of $X$: since $W$ is quasi-isometric to its index $6$ subgroup $W_0=\Z^2$, we deduce that $Y$ is biLipschitz to $\R^2$. According to Corollary~\ref{cor:asymptotic_cone_injective}, we know that $Y$ is injective. Since $W_0$ is abelian, we know that the asymptotic cone $\R^2$ of $W_0$ acts by isometries on $Y$, so $Y$ is a $2$-dimensional normed vector space. The only such injective normed vector space is isometric to $(\R^2,\ell^\infty)$. However, the group $\frak{S}_3$ acts faithfully by linear isometries on $Y$, but $\frak{S}_3$ does not embed into the linear isometry group $\Z/4\Z \rtimes \Z/2\Z$ of $(\R^2,\ell^\infty)$. This is a contradiction.
\ep

\section{Bicombings and normal forms} \label{sec:bicombings_normalforms}

In general, injective metric spaces and Helly graphs are not uniquely geodesic. It therefore appears very useful to have nice choices of geodesic, or preferred paths, between any pair of points.

\subsection{Geodesic bicombings}

A central tool in the study of injective metric spaces is the notion of geodesic bicombing. Indeed, archetypical examples of injective metric spaces are normed vector spaces with the $\ell^\infty$ norm, which are not uniquely geodesic. However, for these examples, there is a "best" choice of geodesic between any pair of points, namely the affine geodesic. The notion of geodesic bicombings aims at solving the problem of non-uniqueness of geodesics by choosing good geodesics.

\bdf
Let $X$ denote a geodesic metric space. A (geodesic) \emph{bicombing} on $X$ is a choice of geodesics, i.e. a map $\sigma:X \times X \times [0,1] \ra X$ such that, for each $x,y \in X$, the map $t \in [0,1] \mapsto \sigma(x,y,t)$ is a constant speed geodesic from $\sigma(x,y,0)=x$ to $\sigma(x,y,1)=y$. The bicombing $\sigma$ is called:
\bit
\item \emph{reversible} if $\forall x,y \in X, \forall t \in [0,1], \sigma(y,x,1-t)=\sigma(x,y,t)$.
\item \emph{conical} if $\forall x,y,x',y' \in X, \forall t \in [0,1]$, we have
$$d(\sigma(x,y,t),\sigma(x',y',t)) \leq (1-t)d(x,x')+td(y,y').$$
\item \emph{convex} if $\forall x,y,x',y' \in X$, the map $t\in [0,1] \mapsto d(\sigma(x,y,t),\sigma(x',y',t))$ is convex.
\item \emph{consistent} if $\forall x,y \in X, \forall\, 0 \leq a \leq b \leq 1, \forall t \in [0,1]$, we have
$$\sigma(x,y,(1-t)a+tb)=\sigma(\sigma(x,y,a),\sigma(x,y,b),t).$$
\eit
\edf

\brk
Note that any convex bicombing is conical. Moreover, any conical, consistent bicombing is also convex. However, there exist convex bicombings which are not consistent, see for instance~\cite{basso_miesch_convex_non_consistent_bicombings}.
\erk

\bthm \cite[Proposition~3.8]{lang}
Any injective metric space admits a reversible conical geodesic bicombing, which is equivariant by isometries.
\ethm

\bp
According to Lemma~\ref{lem:delta_retract}, if $X$ is an injective metric space, then there exists a $1$-Lipschitz retraction $p:\Delta(X) \ra X$ that is equivariant under $\Isom(X)$. Now, for any $x,y \in X$ and $t \in [0,1]$, let us define $\sigma(x,y,t)=\pi((1-t)e(x)+te(y))$, where $e:X \ra \Delta(X)$ is the standard embedding $x \mapsto d(x,\cdot)$. It is clear that $\sigma$ is a reversible and conical bicombing.
\ep

\brk
Any metric space with a conical geodesic bicombing is contractible. This provides another point of view on the contractibility of injective metric spaces, proved in Corollary~\ref{cor:injective_contractible}.
\erk

In order to state the following result, we need to introduce the notion of combinatorial dimension of a metric space.

\bdf
Let $X$ denote a metric space. Its \emph{combinatorial dimension} is the topological dimension of its injective hull.
\edf

Note that, in practice, most examples of injective metric spaces with finite topological dimension will be cell complexes, for which the topological dimension coincides with the maximal dimension of cells.

\bthm[Descombes-Lang, Bicombings, Theorem~4.1]
A metric space $X$ has combinatorial dimension at most $n$ if and only if, for set $Z \subset X$ with cardinality $2(n+1)$, and for any fixed point free involution $i:Z \ra Z$, there exists a fixed point free bijection $j:Z \ra Z$ distinct from $i$ such that:
$$\sum_{z \in Z} d(z,i(z)) \leq \sum_{z \in Z} d(z,j(z)).$$
\ethm

\bexe (Exercise)
A metric space is $0$-hyperbolic if and only if it has combinatorial dimension $1$.
\eexe

\bq
Is there a local characterization of combinatorial dimension?
\eq

\bq
What is the combinatorial dimension of a median graph (i.e. the $1$-skeleton of a CAT(0) cube complex)?
\eq

\bq
What is the combinatorial dimension of Euclidean buildings?
\eq

We can now cite the most precise statements about convex bicombings for injective metric spaces, due to Descombes and Lang.

\bthm \cite[Theorem~1.2]{descombes_lang_hyperbolicity} \label{thm:atmost_one_bicombing}
Let $X$ denote a metric space with finite combinatorial dimension. Then $X$ admits at most one convex geodesic bicombing, which is then consistent and reversible.
\ethm

\bthm \cite[Theorem~1.1]{descombes_lang_hyperbolicity} \label{thm:exists_one_bicombing}
Any proper injective metric space admits at least one consistent geodesic bicombing.
\ethm

\bcor
Let $X$ denote a proper injective metric space with finite combinatorial dimension. Then $X$ admits a unique convex consistent reversible geodesic bicombing.
\ecor

\bexe (Exercise)
Let us consider the union $X$ of a square and a $3$-cube along an edge (see Figure~\ref{fig:geodesics_cube_lp}). Let us endow $X$ with the piecewise $\ell^p$ metric, for $p \in [1,\infty]$. Show that, for $p \in (1,\infty)$, the unique geodesic between the opposite vertices $a$ and $b$ intersects the edge $e$ in a different point.

Note that the piecewise $\ell^p$ metrics on CAT(0) cube complexes have been studied in~\cite{haettel_hoda_petyt_Lp_metrics}.
\eexe

\begin{figure}[H]
\begin{center}
\begin{tikzpicture}
\def \p {0.05}
\def \t {(0.5,1)}
\def \op {0.5}
\def \gris {black!10}
\draw[fill] (0,0) circle (\p) node(a) {};
\draw[fill] (2,0) circle (\p) node(b) {};
\draw[fill] (4,0) circle (\p) node(c) {};
\draw[fill] (4,2) circle (\p) node(d) {};
\draw[fill] (2,2) circle (\p) node(e) {};
\draw[fill] (0,2) circle (\p) node(f) {};
\draw[fill] (2,0)+\t circle (\p) node(b') {};
\draw[fill] (4,0)+\t circle (\p) node(c') {};
\draw[fill] (4,2)+\t circle (\p) node(d') {};
\draw[fill] (2,2)+\t circle (\p) node(e') {};

\draw (a.center) -- (b.center) -- (c.center) -- (d.center) -- (e.center) -- (f.center) -- (a.center);
\draw (b.center) -- (e.center);
\draw (b'.center) -- (c'.center) -- (d'.center) -- (e'.center) -- (b'.center);
\draw (b.center) -- (b'.center);
\draw (c.center) -- (c'.center);
\draw (d.center) -- (d'.center);
\draw (e.center) -- (e'.center);
\node at (-0.2,-0.2) {\bfseries $a$};
\node at (4.7,3.2) {\bfseries $b$};
\draw[thick,blue] (a.center) -- (2,0.7) -- (d'.center);

\end{tikzpicture}
\end{center}
\caption{Geodesics in a cube complex with the $\ell^p$ metric}
\label{fig:geodesics_cube_lp}
\end{figure}

Note that the existence of bicombings is useful to define visual boundaries (see~\cite{descombes_lang_hyperbolicity}) and also barycenters, as we will see in Section~\ref{sec:circumcenters_barycenters}.

\subsection{Normal forms for Helly graphs}

In a Helly graph $X$, there are at least two notions of normal forms, i.e. given any two vertices of $X$, a natural choice of a sequence of cliques joining the two vertices.

\mk

The first one, described in~\cite{helly_groups}, satisfies a fellow travelling property. This results in a biautomatic structure on $X$, in case $X$ is locally finite. The important consequence is that Helly groups are biautomatic (\cite[Theorem~1.5]{helly_groups}).

\mk

The second one, described in~\cite{haettel_osajda_locally_elliptic}, does not satisfy the fellow travelling property but is more flexible and still enjoy a local-to-global property (\cite[Theorem~S]{haettel_osajda_locally_elliptic}). This allows to study locally elliptic groups of automorphisms of Helly graphs (\cite[Theorem~H]{haettel_osajda_locally_elliptic}).

\section{Subdivisions of Helly graphs} \label{sec:subdivisions_helly}

In a CAT(0) cube complex, the combinatorics of isometries become simpler when passing to the cubical barycentric subdivision, see~\cite{haglund}. We will describe a similar process for Helly graphs. In the process, we will describe a very explicit simplicial structure on the injective hull of a Helly graph, following Lang's work (see~\cite{haettel_helly_automorphisms_subdivisions}).

\mk

\bdf[$N$th Helly subdivision]
Let $X$ denote a Helly graph. For each $N \geq 1$, we will consider the graph $X'_N$ with vertex set
$${X'_N}^{(0)} = EX \cap \left(\f{1}{2N!}\Z\right)^X,$$
with an edge between $f,g \in {X'_n}^{(0)}$ if and only if $d(f,g) = \f{1}{2N!}$. The graph $X'_N$ will be called the \emph{$N^\text{th}$ Helly subdivision} of $X$. When $N=1$, we will also denote $X'_1=X'$ the first Helly subdivision. \edf

Note that the factor $\f{1}{2N!}$ is present to ensure a nesting property of the various Helly subdivisions, and also because it provides a more natural statement for the classification of automorphisms of Helly graphs, see Theorem~\ref{thm:dichotomoy_automorphisms_helly}

\bthm
For any Helly graph $X$ and any $N \geq 1$, the $N^\text{th}$ Helly subdivision $X'_N$ is a Helly graph, and the inclusion $X \ra X'_N$ is a $2N!$-homothetic embedding.
\ethm

\bp
This is an immediate consequence of~\cite[Theorem~4.4]{helly_groups}, applied to the discrete metric space $(X,2N!d_X)$.
\ep

\mk

The $N^\text{th}$ Helly subdivision will mostly be useful for the study of automorphisms of Helly graphs, in Section~\ref{sec:classification_automorphisms}. However, the first Helly subdivision is much more natural, and has equivalent descriptions, which we will describe now. When studying actions by automorphisms on Helly graphs, it looks often useful to pass to the first Helly subdivision.

\mk

\bdf[Round cliques]
If $X$ is a graph, we say that a clique $\sigma \subset X$ is \emph{round} if it is an intersection of balls of $X$.
\edf

\bthm \label{thm:characterizations_first_helly_subdivision}
Let $X$ denote a Helly graph, and let $X'$ denote the first Helly subdivision of $X$. Then the vertex set of $X'$ is
\beq X'^{(0)} &=& E(X) \cap \left( \f12 \N\right)^X\\
&=& \{\mbox{round cliques of }X\}\\
&=& X \cup \{\mbox{non-empty intersections of maximal cliques of }X\}.\eeq
Two vertices $\sigma,\tau$ of $X'$, considered as subsets of $X$ in either description, are adjacent in $X'$ if and only if $\sigma \cap \tau \neq \emptyset$ and $\sigma \cup \tau$ is a clique of $X$.
\ethm

The main technical point in the proof of the theorem is the following lemma.

\blem \label{lem:identification_round_cliques}
Let $X$ denote a Helly graph, let $X'$ denote the first Helly subdivision, and let $P_X$ denote the set of round cliques of $X$. The following map is a bijection:
\beq \sigma:X' & \ra & P_X \\
p & \mapsto & \sigma(p)=\bigcap_{x \in X} B(x,\lceil d(p,x) \rceil).\eeq
\elem

For the proofs of Theorem~\ref{thm:characterizations_first_helly_subdivision} and Lemma~\ref{lem:identification_round_cliques}, we refer the reader to~\cite{haettel_helly_automorphisms_subdivisions}.

\mk

We are now able to give a more precise description of the relationship between the Helly hull and the injective hull of a connected graph. Note that Lang has a precise and subtle description of the injective hull of any graph (see~\cite{lang}).

\bpro \label{pro:injective_hull_of_helly_hull}
Let $X$ denote an arbitrary connected graph, and let $H(X)$ denote the Helly hull of $X$. Then the inclusion $X \ra H(X)$ uniquely extends to an isometry between the injective hulls of (the vertex sets of) $X$ and of $H(X)$.
\epro

\bp
It suffices to remark that the vertex set of $H(X)$ may be defined as $E(X) \cap \Z^X$, so in particular $E(X)$ is an injective metric space containing isometrically $H(X)$. By the minimality property of the injective hull, we conclude that $E(H(X))=E(X)$.
\ep

In the case where the graph is already a Helly graph, there is a very simple simplicial structure on the injective hull, that appears as a simple refinement of the cellular structure (see~\cite{dress,lang,dress_moulton_terhalle_Ttheory,chrobak_larmore_injective_hull_algorithm}).

\bthm \label{thm:simplicial_structure_injective_hull_helly}
Let $X$ denote a Helly graph with finite combinatorial dimension. Then the injective hull $E(X)$ has a simplicial complex structure isomorphic to the topological realization of the poset $P_X$ of all round cliques of $X$, partially ordered by inclusion.
Moreover, $E(X)$ is isometric to the $\ell^\infty$ orthoscheme realization of $P_X$ (scaled by a factor $\f{1}{2}$).
\ethm

We refer to Section~\ref{subsec:orthoscheme} for the definition of the orthoscheme realization.

\bp
According to~\cite[Theorem~4.5]{lang}, the injective hull $EX$ may be realized as an isometric subset of $\R^X$, and the injective hull $EX$ of $X$ has a natural cell decomposition satisfying the following. For each cell $C$ of $EX$, there is a finite set of vertices $x_1,\dots,x_n$ of $X$ such that the map
\beq C & \ra & \R^n \\
p & \mapsto & (d(p,x_1),\dots,d(p,x_n))\eeq
is an isometry (with the $\ell^\infty$ metric on $\R^n$) onto the compact convex subspace of $\R^n$ defined by inequalities of the type
$$ \pm d(\cdot,x_i) \pm d(\cdot,x_j) \leq D,$$
for some $1 \leq i < j \leq n$ and $D \in \Z$, and also of the type
$$\pm d(\cdot,x_i) \leq  D',$$
for some $1 \leq i \leq n$ and $D' \in \f{1}{2}\Z$.
In particular there is an affine structure on $C$. Moreover, for any $x \in X$, for any $p _1,\dots,p_k \in C$ and $t_1,\dots,t_k \geq 0$ such that $t_1+\dots+t_k=1$, we have
$$d(x,\sum_{i=1}^k t_ip_i) = \sum_{i=1}^k t_id(x,p_i).$$

\mk

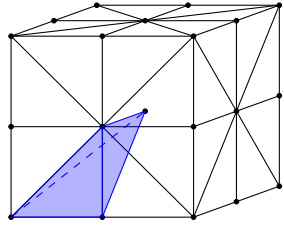
\begin{figure}[H]
\begin{center}
\begin{tikzpicture}[scale = 0.6]
\def \p {0.05}
\def \r {2}
\def \m {1}
\def \a {20}
\def \op {0.5}
\def \gris {black!10}
\draw[fill] (0,0) circle (\p) node(000) {};
\draw[fill] (4,0) circle (\p) node(100) {};
\draw[fill] (0,4) circle (\p) node(010) {};
\draw[fill] (4,4) circle (\p) node(110) {};
\draw[fill] (2,0) circle (\p) node(200) {};
\draw[fill] (0,2) circle (\p) node(020) {};
\draw[fill] (2,4) circle (\p) node(210) {};
\draw[fill] (4,2) circle (\p) node(120) {};
\draw[fill] (2,2) circle (\p) node(220) {};
\draw[fill] (2,2) + (\a:\m) circle (\p) node(222) {};
\draw[fill] (4,0) + (\a:\r) circle (\p) node(101) {};
\draw[fill] (0,4) + (\a:\r) circle (\p) node(011) {};
\draw[fill] (4,4) + (\a:\r) circle (\p) node(111) {};
\draw[fill] (4,0) + (\a:\m) circle (\p) node(102) {};
\draw[fill] (4,2) + (\a:\m) circle (\p) node(122) {};
\draw[fill] (4,2) + (\a:\r) circle (\p) node(121) {};
\draw[fill] (4,4) + (\a:\m) circle (\p) node(112) {};
\draw[fill] (0,4) + (\a:\m) circle (\p) node(012) {};
\draw[fill] (2,4) + (\a:\r) circle (\p) node(211) {};
\draw[fill] (2,4) + (\a:\m) circle (\p) node(212) {};

\draw (000.center) -- (100.center) -- (110.center) -- (010.center) -- (000.center);
\draw (000.center) -- (110.center);
\draw (100.center) -- (010.center);
\draw (200.center) -- (210.center);
\draw (020.center) -- (120.center);

\draw (100.center) -- (101.center) -- (111.center) -- (110.center);
\draw (100.center) -- (111.center) -- (010.center);
\draw (101.center) -- (110.center) -- (011.center);
\draw (010.center) -- (011.center) -- (111.center);
\draw (012.center) -- (112.center) -- (102.center);
\draw (210.center) -- (211.center);
\draw (120.center) -- (121.center);

\draw[blue, dashed] (000.center) -- (222.center);
\draw[blue,fill opacity=0.3,fill=blue] (220.center) -- (000.center) -- (200.center) -- (220.center) -- (222.center) -- (200.center);


\end{tikzpicture}
\end{center}
\caption{The partition of a cube in $\R^3$ into standard orthosimplices.}
\label{fig:partition}
\end{figure}

Note that the hyperplanes of $\R^n$
$$\left\{\pm x_i \pm x_j = D \st 1 \leq i < j, D \in \Z\right\} \mbox{ and } \left\{x_i = D' \st 1 \leq i \leq n, D' \in \f{1}{2}\Z\right\}$$
partition $\R^n$ into (open) standard orthosimplices with edge lengths $\f{1}{2}$, see Figure~\ref{fig:partition}.

\mk

We may consider the refinement of Lang's cell decomposition of $E(X)$, obtained by considering all possible hyperplanes $\{d(\cdot,x) \pm d(\cdot,y) = D\}$, for $x,y \in X$ and $D \in \Z$, and $\{d(\cdot,x) = D'\}$, for $x \in X$ and $D' \in \f{1}{2}\Z$. Each cell from Lang's decomposition is now refined into a finite union of orthoscheme simplices with edge lengths $\f{1}{2}$. Let $O(X)$ denote this simplicial decomposition of $E(X)$.

Notice that every vertex $v \in O(X)$ is such that for each $x \in X$, we have $d(v,x) \in \f12\N$, so $v$ is a vertex of $X'$. Conversely, every vertex $v$ of $X'$ is a vertex of $O(X)$. 

\mk

Moreover, if $v,w$ are vertices of $O(X)$ contained in a common simplex $\sigma$ of $O(X)$, we will prove that the corresponding round cliques $\sigma,\tau \subset X$ are contained in one another. By contradiction, assume that there exists $x \in \sigma \bs \tau$ and $y \in \tau \bs \sigma$. According to the proof of Lemma~\ref{lem:identification_round_cliques}, we deduce that $d(x,\sigma)=\f12$ and $d(x,\tau)=1$, and similarly $d(y,\sigma)=1$ and $d(y,\tau)=\f12$. Hence $\sigma$ and $\tau$ are separated by the hyperplane $\{p \in E(X) \st d(p,x)-d(p,y)=0\}$ of $O(X)$. This contradicts the assumption that $v,w$ are adjacent vertices of $O(X)$.

\mk

Conversely, let us consider two round cliques $\sigma,\tau \subset X$ such that $\sigma \subset \tau$, we will prove that they correspond to adjacent vertices of $O(X)$. It is sufficient to prove that they are not separated by a hyperplane.

Let us fix $x \in X$, $D \in \f12\Z$, since $d(\sigma,\tau)=\f12$, we know that $\sigma$ and $\tau$ are not separated by the hyperplane $\{p \in E(X) \st d(p,x)=D\}$.

Let us fix $x,y \in X$, $\eps=\pm 1$ and $D \in \Z$, and assume by contradiction that  $\sigma$ and $\tau$ are separated by the hyperplane $\{p \in E(X) \st d(p,x)+\eps d(p,y) =D\}$. Since $d(\sigma,\tau)=\f12$, this implies that $d(\sigma,x)+\eps d(\sigma,y) =D \pm \f12$ and $d(\tau,x)+\eps d(\tau,y) =D \mp \f12$. It also implies that $|d(\sigma,x)-d(\tau,x)|=\f12$ and $|d(\sigma,y)-d(\tau,y)|=\f12$. According to the proof of Lemma~\ref{lem:identification_round_cliques}, this implies that there exist $p,q \in \N$ such that, for each $z \in \sigma$, we have $d_X(z,x)=p$ and $d_X(z,y)=q$. Thus $d(\sigma,x)=p$ and $d(\sigma,y)=q$, so $d(\sigma,x)+\eps d(\sigma,y) \neq D \pm \f12$. This is a contradiction.

We conclude that $\sigma$ and $\tau$ are adjacent vertices in $O(X)$.
\ep

In particular, Proposition~\ref{pro:injective_hull_of_helly_hull} and Theorem~\ref{thm:simplicial_structure_injective_hull_helly} provide a simple "recipe" for describing the injective hull of a connected graph $X$:
\ben
\item Determine the Helly hull $H(X)$ of $X$ (see Section~\ref{sec:hulls}).
\item Determine the round cliques of $P_{H(X)}$ of $H(X)$.
\item The injective hull $E(X)$ is the geometric realization of $P_{H(X)}$.
\een

\mk

A nice consequence of the work of Lang, which can also be understood with respect to the above simplicial structure, is the following.

\bcor \label{cor:helly_implies_injective}
Let $G$ denote a Helly group, i.e. a group acting properly and cocompactly by automorphisms on a Helly graph. Then $G$ is an injective group, i.e. $G$ acts properly and cocompactly by isometries on an injective metric space.

Similarly, any coarsely Helly group is coarsely injective.
\ecor

\bp
According to~\cite[Theorem~4.5]{lang}, the injective hull $E(X)$ of a Helly graph $X$ is at distance at most $1$ from $X$, and is a proper metric space.
\ep

We will see in Section~\ref{subsec:orthoscheme_complexes} that the local injectivity of $E(X)$ has a simple interpretation from the poset properties of $P_X$.

\section{Circumcenters, barycenters and fixed points} \label{sec:circumcenters_barycenters}

In CAT(0) spaces, one can naturally define circumcenters of bounded sets and barycenters of probability measures. We can also define such objects in the case of injective spaces and Helly graphs, as explained below.

\bthm[Circumcenter] \label{thm:circumcenter}
Let $X$ denote an injective metric space. Let ${\frak B}(X)$ denote the set of all non-empty bounded subsets of $X$. Then there is a circumcenter map $\gamma: {\frak B}(X) \ra X$ such that:
\ben
\item For each $K \in {\frak B}(X)$ and each ball $B \subset X$ containing $K$, we have $\gamma(K) \in B$.
\item The map $\gamma$ is equivariant with respect to the isometry group of $X$.
\een
\ethm

\bp
Fix $K \in {\frak B}(X)$. Let
$$K_0= \bigcap_{B \mbox{ ball in $X$ containing } K} B,$$
we will define a nested sequence $(K_n)_{n \in \N}$ by induction. Fix $n \in \N$, and let $d_n=\diam(K_n)$ denote the diameter of $K_n$. Let us define
\[K'_{n+1} = \bigcap_{x \in K_n} B\left(x,\frac{d_n}{2}\right) \cap K_n.\]
Note that two balls in the intersection above pairwise intersect, so by the Helly property we deduce that $K'_{n+1}$ is non-empty.

\mk

Recall that the intrinsic radius $ri(A)$ of a subset $A \subset X$ is defined by $ir(A) = \inf\{r \geq 0 \st \exists a \in A, A \subset B(a,r)\}$. Then the intrinsic radius of $K'_{n+1}$ satisfies $ir(K'_{n+1}) \leq \frac{d_n}{2}$. Let us now define
\[ K_{n+1} = \{x \in K'_{n+1} \st K'_{n+1} \subset B(x,ir(K'_{n+1})+2^{-n})\}\]
We have seen that the set $K_{n+1}$ is not empty. Furthermore, we know that the diameter $d_{n+1}$ of $K_{n+1}$ satisfies $d_{n+1} \leq ir(K'_{n+1})+2^{-n} \leq \frac{d_n}{2}+2^{-n}$.

\mk

This defines a nested sequence of $(K_n)_{n \in \N}$ of non-empty bounded subsets of $X$. We also know that, for each $\alpha \in \left(\f{1}{2},1\right)$, we have
\[ \liml_{n \ra \pif} \f{d_n}{\alpha^n} = 0.\]
Furthermore, note that for each $n \in \N$, and for each $x \in K_n$ and $y \in K_{n+1}$ we have $d(x,y) \leq \f{d_n}{2}$. We deduce that any sequence $(x_n)_{n \in \N}$ such that $x_n \in K_n$ for each $n \in \N$ is a Cauchy sequence. Since the injective metric space $X$ is complete, we deduce that any such sequence converges in $X$. In addition, since the diameter of $K_n$ goes to $0$ as $n$ goes to $+\infty$, there exists a unique such limit: let us denote it $\gamma(K)$.

\mk

Fix a ball $B$ in $X$ containing $K$: we have that $K_0 \subset K$. Since the sequence is nested, we have $\gamma(K) \in B$.

\mk

Furthermore, it is clear from the definition of $\gamma$ that it is equivariant with respect to isometries of $X$.
\ep

In injective metric spaces, or more generally in metric spaces with reversible conical bicombings, one can consider the barycenter of a probability measure. These constructions have origins in the works of many authors, including~\cite{essahib_heinich_barycenter,descombes_asymptotic_rank,navas_ergodic_barycentre,petyt_PhD_thesis}.

\bthm[Barycenter] \cite[Theorem~3.4]{basso_fixed_point}
Let $X$ denote an injective metric space, let ${\cal P}(X)$ denote the space of Borel probability measures on $X$, and let $W_1$ denote the Wasserstein distance on ${\cal P}_1(X)$ (see~\cite{basso_fixed_point}). There exists a map $\beta:{\cal P}(X) \ra X$ such that:
\bit
\item For any $x \in X$, we have $\beta(\delta_x)=x$.
\item More generally, for any ball $B \subset X$ and for any measure $\mu \in {\cal P}(X)$ with support in $B$, we have $\beta(\mu) \in B$.
\item The map $\beta$ is $1$-Lipschitz.
\item The map $\beta$ is equivariant with respect to the isometry group of $X$.
\eit
\ethm

\bthm[Fixed points in injective spaces] \cite[Proposition~1.2]{lang}
Let $X$ denote an injective metric space, and let $G$ denote a group acting by isometries on $X$ with bounded orbits. Then $G$ has a fixed point, and the fixed point set $X^G$ is injective.
\ethm

\bp
We will give a different proof than~\cite[Proposition~1.2]{lang}, using circumcenters instead.

Assume that $G$ has bounded orbits, and let $x \in X$. According to Theorem~\ref{thm:circumcenter}, the circumcenter $\gamma(G \cdot x)$ is fixed by $G$.

\mk

Now let us consider a family $(x_i,r_i)_{i \in I}$ in $X^G$, with $\forall i,j \in I, r_i+r_j \geq d(x_i,x_j)$. Since $X$ is injective, the intersection $K=\bigcap_{i \in I} B(x_i,r_i)$ is non-empty. Since it is bounded and $G$-invariant, its circumcenter $\gamma(K)$ is well-defined, and it is fixed by $G$. Moreover, according to Theorem~\ref{thm:circumcenter}, the circumcenter $\gamma(K)$ is contained in any ball containing $K$: we deduce that $\gamma(K) \in \bigcap_{i \in I} B(x_i,r_i) \cap X^G$. Hence $X^G$ is injective.
\ep

We can do similar constructions for Helly graphs.

\bthm[Circumclique] \label{thm:circumclique}
Let $X$ denote a Helly graph. Let ${\frak B}(X)$ denote the set of all non-empty bounded subsets of $X$, and let $P(X)$ denote the set of all cliques of $X$ which are intersections of balls of $X$. Then there is a circumclique map $\gamma: {\frak B}(X) \ra P(X)$ such that:
\ben
\item For each $K \in {\frak B}(X)$ and each ball $B \subset X$ containing $K$, we have $\gamma(K) \subset B$.
\item The map $\gamma$ is equivariant with respect to the automorphism group of $X$.
\een
\ethm

\bp
Fix $K \in {\frak B}(X)$. Let
$$K_0= \bigcap_{B \mbox{ ball in $X$ containing } K} B,$$
we will define a nested sequence $(K_n)_{n \in \N}$ by induction. Fix $n \in \N$, and let $d_n=\diam(K_n)$ denote the diameter of $K_n$. Let us define
\[K'_{n+1} = \bigcap_{x \in K_n} B\left(x,\left\lceil \frac{d_n}{2}\right\rceil \right) \cap K_n.\]
Note that two balls in the intersection above pairwise intersect, so by the Helly property we deduce that $K'_{n+1}$ is non-empty.

\mk

Recall that the intrinsic radius $r(A)$ of a subset $A \subset X$ is defined by $ir(A) = \min\{r \in \N \st \exists a \in A, A \subset B(a,r)\}$. Then the intrinsic radius of $K'_{n+1}$ satisfies $ir(K'_{n+1}) \leq \left\lceil \frac{d_n}{2}\right\rceil$. Let us now define
\[ K_{n+1} = \{x \in K'_{n+1} \st K'_{n+1} \subset B(x,ir(K'_{n+1}))\}\]
Note that, since $K'_{n+1}$ is an intersection of balls in the injective metric space $X$, the set $K_{n+1}$ is not empty. Furthermore, we know that the diameter $d_{n+1}$ of $K_{n+1}$ satisfies $d_{n+1} \leq ir(K'_{n+1}) \leq \left\lceil \frac{d_n}{2}\right\rceil$.

\mk

This defines a nested sequence of $(K_n)_{n \in \N}$ of non-empty bounded subsets of $X$. So there exists $n \in \N$ such that for every $k \geq n$, we have $d_k=1$, and thus $K_k=K_n$. Let us define $\beta(K)=K_n$: it is a clique of $X$, which is an intersection of balls of $X$, so $\beta(K) \in P(K)$.

\mk

Furthermore, it is clear from the definition of $\gamma$ that it is equivariant with respect to isometries of $X$.
\ep

\bthm[Fixed points in Helly graphs]
Let $X$ denote a Helly graph, and let $G$ denote a group acting by automorphisms on $X$ with bounded orbits, and let $X'$ denote the first Helly subdivision of $X$. Then $G$ stabilizes a round clique of $X$, $G$ has a fixed vertex in $X'$, and $X'^G$ is a Helly graph.
\ethm

\bp
Assume that $G$ has bounded orbits, and let $x \in X$. According to Theorem~\ref{thm:circumclique}, the circumclique $\gamma(G \cdot x)$ is stabilized by $G$. According to Theorem~\ref{thm:characterizations_first_helly_subdivision}, the round clique $\gamma(G \cdot x)$ corresponds to a vertex of the first Helly subdivision $X'$, so $G$ has a fixed point in $X'$.

\mk

Assume that $(\sigma_i)_{i \in I}$ are pairwise intersecting maximal cliques of $X'^G$. Then each $\sigma_i$ is a family of pairwise intersecting $G$-invariant cliques of $X$ (which are intersections of balls). Let us define $\tau_i = \cap \sigma_i$: it is a clique of $X$, which is an intersection of balls of $X$, so $\tau_i \in P(X)$, and it is also $G$-invariant. Moreover, since $\sigma_i$ is maximal, we know that $\tau_i \in \sigma_i$. Moreover, for each $i,j \in J$, since $\sigma_i$ and $\sigma_j$ do intersect in $X'$, we deduce that $\tau_i$ and $\tau_j$ intersect in $X$.

So we have a family of pairwise intersecting $G$-invariant cliques $(\tau_i)_{i \in I}$ of $X$, which are intersections of balls: since $X$ is Helly, the intersection $\tau = \bigcap_{i \in I} \tau_i$ is non-empty and $G$-invariant. The clique $\tau$ of $X$ is an intersection of balls of $X$, so $\tau \in P(X)$. Moreover, for each $i \in I$, we have $\tau \subset \tau_i$, so $\tau$ is adjacent to $\sigma_i$ in $X'$. This implies that the vertex $\tau \in X'^G$ is contained in the intersection $\bigcap_{i \in I} \sigma_i$. In other words, $X'^G$ is clique-Helly.

\mk

Let us now prove that $X'^G$ is connected and simply connected.

\mk

More generally, assume that a group $G$ acts on a Helly graph $X$ such that the following hold:
\bit
\item $G$ fixes a vertex of $X$.
\item The fixed point set $X^G$ is clique-Helly.
\item Stabilizers of cliques of $X$ are pointwise stabilizers.
\eit
We will prove that $X^G$ is connected and simply connected. We will start by proving that $X^G$ is isometrically embedded in $X$. Let us fix vertices $x,y \in X^G$ at distance $d_X(x,y)=n$ in $X$: we will prove, by induction on $n$, that $d_{X^G}(x,y)=n$.

For $n=1$ it is clear, so let us assume that $n \geq 2$. Let $A$ denote the set of all vertices $z \in X$ such that $d_X(x,z)=n-1$ and $d_X(z,y)=1$: it is a non-empty, bounded, $G$-invariant set. So its circumclique $\gamma(A)$ is such that, for any $z \in A$, we have $d_X(x,z)=n-1$ and $d_X(z,y)=1$. Moreover, since $A$ is $G$-invariant and by assumption on $X$, we deduce that $z$ is fixed by $G$. Hence by induction we have $d_{X^G}(x,z)=n-1$ and $d_{X^G}(x,y)=n$.

Hence $X^G$ is connected.

\mk

Let us now 

There exists a unique clique $\sigma \subset X$, maximal with the following property, such that $d(x,\sigma)=n-1$ and $d(y,\sigma)=1$.

vertex $x$ is also 

 adjacent to all $\tau_i$, for $i \in I$. In particular, $x$ is adjacent to every clique in every $\sigma_i$, for $i \in I$

Let us consider a family $(\sigma_i,r_i)_{i \in I}$ in $X'^G$, with $r_i \in \N$, with $\forall i,j \in I, r_i+r_j \geq d_{X'}(\sigma_i,\sigma_j)$. For each $i \in I$, we have that $\sigma_i$ is a clique of $X$ which is an intersection of balls of $X$. Moreover, for each $i,j \in I$, we have 

Since $X$ is Helly, the intersection $K=\bigcap_{i \in I} B(x_i,r_i)$ is non-empty. Since it is bounded and $G$-invariant, its circumcenter $\gamma(K)$ is well-defined, and it is fixed by $G$. Moreover, according to Theorem~\ref{thm:circumcenter}, the circumcenter $\gamma(K)$ is contained in any ball containing $K$: we deduce that $\gamma(K) \in \bigcap_{i \in I} B(x_i,r_i) \cap X^G$. Hence $X^G$ is injective.
\ep

\section{Classification of automorphisms} \label{sec:classification_automorphisms}

In order to study isometries of injective metric space of automorphisms of Helly graphs, it is useful to have a precise description of their behaviour.

\mk

Descombes and Lang proved general results for isometries of metric spaces with a convex bicombing (see~\cite{descombes_lang_flats}).

\bthm[Isometries of injective metric spaces] \cite[Proposition~5.5]{descombes_lang_flats}
Let $X$ denote a proper injective metric space, and let $G$ denote a group acting properly and cocompactly on $X$. Then, for any element $g \in G$, we have the following dichotomy:
\bit
\item Either $g$ is elliptic, i.e. $g$ satisfies one of the following equivalent properties:
\bit
\item $g$ has bounded orbits in $X$.
\item $g$ fixes a point in $X$.
\eit
\item Or $g$ is hyperbolic, i.e. $g$ satisfies one of the following equivalent properties:
\bit
\item For some point $x \in X$, the map $n \in \Z \mapsto g^n \cdot x$ is a quasi-isometric embedding.
\item For some (equivalently, any) $G$-invariant convex consistent geodesic bicombing $\sigma$ on $X$, the isometry $g$ translates a $\sigma$-geodesic axis in $X$.
\eit
\eit
\ethm

\bp
This is essentially stated in~\cite[Proposition~5.5]{descombes_lang_flats}, and the existence of a $G$-invariant convex geodesic bicombing follows from~\cite[Theorem~1.1]{descombes_lang_hyperbolicity} (also stated as Theorem~\ref{thm:exists_one_bicombing}).
\ep

\brk
Parabolic isometries of injective metric spaces are abundant, simply because every metric space embeds into its injective hull. For instance, let us consider an action of $\SL(2,\Z)$ by isometries on $\H^2$ properly with finite covolume. Then $\SL(2,\Z)$ acts properly by isometries on the injective hull $X$ of $\H^2$ (though not cocompactly), and unipotent elements of $\SL(2,\Z)$ act as parabolic isometries on $\H^2$ and on $X$.
\erk

We now turn to the study of automorphisms of Helly graphs. Similarly to the injective case, we need an assumption to forbid parabolic automorphisms, which is provided by the combinatorial dimension.

\bthm \label{thm:dichotomoy_automorphisms_helly}
Let $\Gamma$ denote a Helly graph of combinatorial dimension $N$, and let $g$ denote an automorphism of $\Gamma$. Then we have the following dichotomy:
\bit
\item Either $g$ is elliptic, i.e. $g$ satisfies one of the following equivalent properties:
\bit
\item $g$ has bounded orbits in $\Gamma$.
\item $g$ stabilizes a clique in $\Gamma$. 
\item $g$ fixes a point in the injective hull $E\Gamma$.
\item $g$ fixes a vertex of the first Helly subdivision $\Gamma'$ of $\Gamma$.
\eit
\item Or $g$ is hyperbolic, i.e. $g$ satisfies one of the following equivalent properties:
\bit
\item For some vertex $v \in \Gamma$, the map $n \in \Z \mapsto g^n \cdot v$ is a quasi-isometric embedding.
\item $g$ translates a geodesic axis in the injective hull $E\Gamma$.
\item There exists a vertex $v$ of the first Helly subdivision $\Gamma'$ of $\Gamma$, $1 \leq a \leq 2N$ and $L \in \N \bs \{0\}$ such that, for any $n \in \Z$, we have $d(g^{an} \cdot v,v) = |n|L$. 
\item There exists a vertex $v$ of the $N$th Helly subdivision $\Gamma'_N$ of $\Gamma$ and $L \in \N \bs \{0\}$ such that, for any $n \in \Z$, we have $d(g^n \cdot v,v) = |n|L$. 
\eit
\eit
\ethm

\bp
See~\cite{haettel_helly_automorphisms_subdivisions}.
\ep

\brk
Note that the restriction on the combinatorial dimension is necessary to have the dichotomy. For instance, consider $g$ is a parabolic automorphism of a connected graph $\Gamma$: to be more specific, one may consider a finitely generated group $G$, $\Gamma$ a Cayley graph of $G$, and $g$ an infinite order element of $G$ which is distorted. Then $g$ extends to an automorphism of the Helly hull of $\Gamma$, which is neither elliptic nor loxodromic.
\erk

\brk
This result is very similar to the ananlogous statement for CAT(0) cube complexes by Haglund (see~\cite{haglund}), which is stated in terms of the cubical barycentric subdivision.
\erk

\brk
Note that it is not sufficient to consider only the first Helly subdivision for hyperbolic automorphisms: consider the Helly graph $\Gamma$ with vertex set $\Z^N$, with the standard Helly structure. Let $g$ denote the following automorphism of $\Gamma$:
$$g \cdot (x_1,x_2,\dots,x_N) = (x_2+1,x_3,x_4,\dots,x_N,x_1).$$
Then $g$ is hyperbolic with translation length $\f{1}{N}$, whereas any automorphism with combinatorial axis in the first Helly subdivision has translation length in $\f{1}{2}\N$.
\erk

\section{Constructions of injective spaces and Helly graphs} \label{sec:constructions}

We will now describe some simple methods to build injective metric spaces and Helly graphs, notably from median spaces and from thickening of cell complexes. We will also discuss why hyperbolic groups are Helly.

\subsection{From median to injective metric spaces}

Bowditch and Miesch (independently) give a construction of an injective metric for many metric median spaces. Let us recall the following definition (see for instance~\cite{bowditch_median_metric,bowditch_book_median}).

\bdf[Metric median space]
A metric space $X$ is \emph{metric median} if, for any $x,y,z \in X$, there exists a unique $m=\mu(x,y,z)$, called the \emph{median} of $\{x,y,z\}$, such that
\beq d(x,y)&=&d(x,m)+d(m,y)\\
d(y,z)&=&d(y,m)+d(m,z)\\
d(x,z)&=&d(x,m)+d(m,z).\eeq
Then \emph{rank} of $X$ is the supremum of all $r \in \N$ such that there exists an homothetic embedding of the $r$-cube $(\{0,1\}^r,\ell^1)$ into $X$.
\edf 

\bexes\

\bit
\item The real line $\R$ is metric median, with median $\mu_\R(x,y,z)=y$ if $x \leq y \leq z$.
\item For each $n \in \N$, the vector space $(\R^n,\ell^1)$ is metric median of rank $n$.
\item For each measured space $\Omega$, the vector space $L^1(\Omega)$ is metric median (usually of infinite rank).
\item Any $n$-dimensional CAT(0) cube complex $X$, endowed with the piecewise $\ell^1$ metric, is metric median of rank $n$.
\eit
\eexes

One can easily show that, in a metric median space, compact convex subsets satisfy the Helly property (see~\cite{bowditch_median_injective}).

\bthm \label{thm:convex_subsets_helly}
Let $X$ denote a complete, connected, finite rank metric median space, and let $\mu:X^3 \ra X$ denote the median of $X$. Then the family of $\mu$-convex closed bounded subsets of $X$ satisfies the Helly property.
\ethm

\bp
We will only give the proof for $\mu$-convex compact subsets of $X$. We will use transfinite induction on the ordinal $I$ to prove that any family of pairwise intersecting convex compact subsets $(C_i)_{i \in I}$ of $X$ has a non-empty intersection.

\mk

Let us start with $|I|=3$: let $C_1,C_2,C_3$ denote pairwise intersecting compact convex subsets of $X$. For each $1 \leq i \leq 3$, let $x_i \in C_i \cap C_{i+1}$ (modulo $3$), and let $y=\mu(x_1,x_2,x_3)$: by convexity, we have $y \in C_1 \cap C_2 \cap C_3$.

\mk

Assume now that $I$ is an arbitrary ordinal for which $X$ is $I$-hyperconvex, and consider it successor $I'=I+1$. Consider a family of pairwise intersecting compact convex subsets $(C_i)_{i \in I'}$. For each $i \in I$, let us define $C'_i=C_i \cap C_{I+1}$: it is a compact convex subset of $X$. According to the $3$-hyperconvexity of $X$, the subsets $(C'_i)_{i \in I}$ pairwise intersect. Hence the convex subsets $(C_i)_{i \in I'}$ have a global intersection.

\mk

Assume now that $I$ a limit ordinal, and consider a family of pairwise intersecting convex compact subsets $(C_i)_{i \in I}$ of $X$. By induction, for each $i \in I$, the intersection $\bigcap_{j \leq i} C_j$ is compact, convex and non-empty. By compactness, each $\bigcap_{j \leq i} C_j$ is compact, so the intersection $\bigcap_{i \in I} C_i$ is non-empty.
\ep

One nice consequence is that the thickening of the $1$-skeleton of a locally finite CAT(0) cube complex is Helly. We will see in the sequel other ways to prove this result, notably to remove the local finiteness assumption.

\bcor
Let $X$ denote a locally finite CAT(0) cube complex. Let $Y$ denote its \emph{thickening}, i.e. the graph with vertex set $X^{(0)}$ with an edge between two vertices if they are contained in a common cube of $X$. Then $Y$ is a Helly graph.
\ecor

\bp
This is a consequence of the fact that balls in $Y$ are finite and median-convex. It is not too hard to prove, and we admit this result here.
\ep

In order to define injective metrics on metric median spaces, Bowditch introduces the notion of contraction.

\bdf[Contraction]
Let $(X,d,\mu)$ denote a metric median space. A \emph{contraction} is a $1$-Lipschitz, median map $f:X \ra \R$, i.e. such that
$$\forall x,y,z \in X, \mu_\R(f(x),f(y),f(z)) = f(\mu(x,y,z)).$$
\edf

\bthm[\cite{bowditch_median_injective}]
Let $(X,d_1)$ denote a connected, complete metric median space of finite rank. Let us define the following quantity, for each $x,y \in X$:
$$d_\infty(x,y) = \sup_{f:X \ra \R \emph{ contraction}} |f(x)-f(y)|.$$
Then $d_\infty$ is an injective metric on $X$, that is biLipschitz to $d_1$. 
\ethm

\bp
It is easy to see that $d_\infty$ is a pseudometric. We will only sketch the argument that $d_\infty$ is injective.

\mk

From the definition, it is easy to see that balls for $d_\infty$ are median-convex: let $x \in X$, $r>0$, and $y,z \in B_{d_\infty}(x,r)$. Let $m=\mu(x,y,z) \in X$, we want to prove that $d_\infty(m,x) \leq r$. For each contraction $f:X \ra \R$, we have $|f(x)-f(y)| \leq r$ and $|f(x)-f(z)| \leq r$, hence we deduce that $|f(x)-f(m)| \leq r$. In conclusion, $m \in B_{d_\infty}(x,r)$.

According to Theorem~\ref{thm:convex_subsets_helly}, we conclude that balls in $(X,d_\infty)$ satisfy the Helly property.

\mk

The hard part is to prove that $d_\infty$ separates points, and that $d_\infty$ is geodesic. This relies on the fine properties of intervals in finite rank metric median spaces.
\ep

This result applies in particular to CAT(0) cube complexes.

\bcor[\cite{bowditch_median_injective},\cite{miesch_CCC}]
Let $X$ denote a locally finite-dimensional CAT(0) cube complex, endowed with the piecewise $\ell^\infty$ metric. Then $X$ is injective.
\ecor

In the sequel, we will see another way to prove this result, using orthoscheme realizations of semilattices (Corollary~\ref{cor:CCC_injective_semilattices}), and also another way to prove that the corresponding graph is Helly, using cell-Helly complexes (Theorem~\ref{thm:CCC_Helly_cell_Helly}).

Using a similar idea, one can prove the following.

\bthm[\cite{haettel_hoda_petyt}]
The mapping class group of a finite type surface, or more generally any hierarchically hyperbolic group, has a proper and cobounded action on an injective metric space.
\ethm

\subsection{Thickening of complexes}

\bdf
Let $X$ denote a graph, with a family $(X_i)_{i \in I}$ of subsets of $X$, called \emph{generalized cells} of $X$. The \emph{thickening} of $X$ (with respect to this family) is the graph with vertex set $X^{(0)}$, with an edge between $u,v \in X^{(0)}$ if they are contained in a common generalized cell of $X$.
\edf

\bthm \cite[Theorem~3.9]{helly_groups} \label{thm:thickening_helly}
Let $X$ denote a combinatorial complex (see~\cite[Chapter~I.8.A]{bridson_haefliger}), and let $(X_i)_{i \in I}$ denote a family of finite full combinatorial subcomplexes of $X$ covering $X$, called generalized cells. Assume that the following hold:
\bit
\item $X$ is simply connected.
\item Each non-empty intersection of generalized cells is connected and simply connected.
\item Generalized cells satisfy the finite Helly property.
\item {\bf Flag condition} For any pairwise intersecting generalized cells $X_1,X_2,X_3$, there exists $i \in I$ such that
$$(X_1 \cap X_2) \cup (X_2 \cap X_3) \cup (X_3 \cap X_1) \subset X_i.$$
\item {\bf Locally bounded} There is no infinite strictly increasing chain of generalized cells.
\eit
The complex $X$ is called \emph{cell Helly}. Then the thickening of $X$ is a Helly graph, and its maximal cliques are generalized cells.
\ethm

Here is a very simple application of this construction to cube complexes.

\bthm \label{thm:CCC_Helly_cell_Helly}
Let $X$ denote a cube complex without infinite cubes. Then $X$ is CAT(0) if and only if its thickening is a Helly graph.
\ethm

\bp
Let us consider all cubes as generalized cells. Then Gromov's flag condition is equivalent to the flag condition from Theorem~\ref{thm:thickening_helly}.
\ep

Another important application of this construction concerns spherical type (or FC type) Artin groups, and more generally Garside groups: see Section~\ref{subsec:garside}.

\subsection{Injective hulls of hyperbolic groups}

We will now give a brief of account of Lang's result that Gromov-hyperbolic groups are Helly.

Recall that a possibly non-geodesic metric space $X$ is called $\delta$-hyperbolic if
$$\forall x,y,z,t \in X, d(x,y)+d(z,t) \leq \max\left(d(x,z)+d(y,t),d(x,t)+d(y,z)\right)+\delta.$$

\bpro \cite[Proposition~1.3]{lang} \label{pro:hull_hyperbolic_space}
Let $X$ denote a $\delta$-hyperbolic metric space, then its injective hull $E(X)$ is $\delta$-hyperbolic as well. If in addition $X$ is geodesic or discretely geodesic, then the image of $X$ in its injective hull is cobounded.
\epro

\bp
Let us consider $e,f,g,h \in E(X)$, and fix $\eps>0$. Let $x \in X$ such that $\|e-f\|_\infty \leq |e(x)-f(x)|+\eps$, so up to exchanging the roles of $e$ and $f$ assume that $\|e-f\|_\infty \leq e(x)-f(x)+\eps$. Let $y \in X$ such that $e(x)+e(y) \leq d(x,y)+\eps$. So
$$\|e-f\|_\infty \leq d(x,y)-f(x)-e(y)+2\eps.$$
Similarly, there exist $z,t \in X$ such that
$$\|g-h\|_\infty \leq d(z,t)-h(z)-g(t)+2\eps.$$
Let us denote $\Sigma=f(x)+e(y)+h(z)+g(t)$. So we deduce that
\beq \|e-f\|_\infty + \|g-h\|_\infty & \leq & d(x,y)+d(z,t)-\Sigma+4\eps\\
& \leq &  \max\left(d(x,z)+d(y,t),d(x,t)+d(y,z)\right)-\Sigma+\delta+4\eps.\eeq
Note that
\beq d(x,z)+d(y,t) - \Sigma &=& (d(x,z)-f(x)-h(z)) + (d(y,t)-e(y)-g(t))\\
&\leq& (f(z)-h(z))+(e(t)-g(t)) \leq \|f-h\|_\infty+\|e-g\|_\infty,\eeq
and similarly $d(x,t)+d(y,z) - \Sigma \leq \|f-g\|_\infty+\|e-h\|_\infty$. So we deduce that
$$\|e-f\|_\infty + \|g-h\|_\infty \leq \max\left(\|f-h\|_\infty+\|e-g\|_\infty,\|f-g\|_\infty+\|e-h\|_\infty\right)+\delta+4\eps.$$
This holds for any $\eps>0$, hence $E(X)$ is $\delta$-hyperbolic.

\mk

Assume now that $X$ is geodesic. Fix $f \in E(X)$, $\eps>0$, and let $x,y \in X$ such that $f(x)+f(y) \leq d(x,y)+\eps$. Since $X$ is geodesic, there exists $z \in [x,y]$ such that $d(x,z) \geq f(x)-\eps$ and $d(z,y) \geq f(y)-\eps$. According to the $\delta$-hyperbolicity of $E(X)$ applied to $\{f,x,y,z\}$, we deduce that
\beq f(z) + d(x,y) & \leq & \max\left(f(x)+d(z,y),f(y)+d(z,x)\right)+\delta\\
& \leq & \max\left(d(x,z)+d(z,y)+\eps,d(y,z)+d(z,x)+\eps \right)+\delta\\
& \leq & d(x,y)+\eps+\delta.\eeq
In particular $f(z) \leq \eps+\delta$, so $f$ is at distance at most $\eps+\delta$ from $X$. The argument when $X$ is discretely geodesic is similar, with different constants.
\ep

A very important consequence concerns Gromov-hyperbolic groups.

\bcor \cite[Theorem~1.4]{lang},\cite[Theorem~1.1]{helly_groups}
Let $G$ denote a Gromov-hyperbolic group, and let $X$ denote a Cayley graph of $G$. Then the Helly hull $H(X)$ is locally finite, and $G$ acts properly cocompactly on $H(X)$. Similarly the injective hull $E(X)$ is proper, finite-dimensional, and $G$ acts properly and cocompactly on $E(X)$. In particular, $G$ is a Helly group.
\ecor

\bp
By $\delta$-hyperbolicity, one deduce that $X$ has stable intervals. According to Theorem~\ref{thm:stable_intervals_proper_hull}, we deduce that $E(X)$ is proper. According to Proposition~\ref{pro:hull_hyperbolic_space}, we know that $E(X)$ is at distance at most $\delta$ from $X$. Hence the action of $G$ on $E(X)$ is proper and cocompact. Moreover, according to Theorem~\ref{thm:exists_helly_hull} we have $H(X)=E(X) \cap \N^X$, hence the action of $G$ on the Helly graph $H(X)$ is also cocompact.
\ep

\subsection{Coarse injectivity, coarse Hellyness} \label{sec:coarse_injectivity}

\bdf[Coarsely injective]
A metric space $X$ is called \emph{coarsely injective} if there exists $\delta>0$ such that for any family of points $(x_i)_{i \in I}$ in X, and for any family of radii $(r_i)_{i \in I}$ in $\R_+$ such that for any $i \neq j$ in $I$, we have $r_i+r_j \geq d(x_i,x_j)$, we require that the intersection $\bigcap_{i \in I} B(x_i,r_i+\delta)$ is non-empty.

A graph with vertex set $X$ is called \emph{coarsely Helly} if there exists $\delta \in \N$ such that for any family of pairwise intersecting combinatorial balls $(B(x_i,r_i))_{i \in I}$, we require that the intersection $\bigcap_{i \in I} B(x_i,r_i+\delta)$ is non-empty.
\edf

Recall that a subset $A$ of a metric space $X$ is \emph{coarsely dense} if there exists $\delta \geq 0$ such that any point of $X$ is at distance at most $\delta$ from a point of $A$.

\bpro \cite[Proposition~3.12]{helly_groups}\ \label{pro:coarsely_injective_coarsely_dense}

Let $X$ denote a metric space. Then the image of $X$ is coarsely dense in its injective hull if and only if $X$ is coarsely injective.

Let $X$ denote the vertex set of a graph. Then the image of $X$ is coarsely dense in its Helly hull if and only if $X$ is coarsely Helly.
\epro

\bp
Assume that the image of $X$ in injective hull $e:X \ra E(X)$ is coarsely dense, i.e. there exists $\delta \geq 0$ such that the $\delta$-neighbourhood of $e(X)$ in $E(X)$ equals $E(X)$. Consider a family of points $(x_i)_{i \in I}$ in X and a family of radii $(r_i)_{i \in I}$ in $\R_+$ such that for any $i \neq j$ in $I$, we have $r_i+r_j \geq d(x_i,x_j)$. Since $E(X)$ is injective and $e$ is an isometric embedding, there exists a point $z \in E(X)$ in the intersection of all balls $B(e(x_i),r_i))$, for $i \in I$. Let $x \in X$ such that $d(e(x),z) \leq \delta$. Then $x \in X$ is in the intersection of all balls $B(x_i,r_i+\delta)$, for $i \in I$. So $X$ is coarsely injective.

\mk

Assume that $X$ is coarsely injective, for a constant $\delta \geq 0$, and fix $z \in E(X)$. For each $x \in X$, let $r_x=d(e(x),z)$. By the triangle inequality, for any $x,x' \in X$, we have $r_x+r_{x'} \geq d(e(x),e(x'))=d(x,x')$. By coarse injectivity, there exists $y \in X$ in the intersection of all balls $B(x,r_x+\delta)$, for $x \in X$. We will prove that $d(z,e(y)) \leq \delta$: assume by contradiction that $d(z,e(y)) > \delta$, let $\eps=\f13(d(z,e(y))-\delta)>0$.

By definition, $d(z,e(y)) = \sup_{x \in X} |d(x,z)-d(x,e(y))|$, so there exists $x \in X$ such that $d(z,e(y)) \leq |d(x,z)-d(x,e(y))|+\eps$. Since $d(x,y) \leq d(x,z)+\delta$, we deduce that $d(z,y) \leq d(x,z)-d(x,y)+\eps$, so $d(x,z) \geq d(x,y)+d(y,z)-\eps$.

According to Theorem~\ref{thm:explicit_description_injective_hull}, there exists $x' \in X$ such that $d(x,z)+d(x',z) \leq d(x,x')+\eps$. Since $d(x',z) \geq d(x',y)-\delta$, we deduce
\beq d(x,x') &\geq& d(x,z)+d(x',z) -\eps \\
& \geq & d(x,y)+d(y,z)+d(x',y)-2\eps-\delta\\
& \geq & d(x,x') +d(y,z)-\delta-2\eps > d(x,x'),\eeq
which is a contradiction.

So the image of $X$ in $E(X)$ is $\delta$-coarsely dense.
\ep

One can use it to deduce the following applicable criterion to decide whether a given group is Helly.

\bcor
Let $G$ denote a group acting properly and cocompactly on a coarsely Helly graph $X$ with stable intervals. Then $G$ acts properly and cocompactly on a Helly graph.
\ecor

\bp
Since $X$ is coarsely Helly, according to Proposition~\ref{pro:coarsely_injective_coarsely_dense}, there exists $\delta \in \N$ such that $X$ is $\delta$-dense in its Helly hull $H(X)$. According to Theorem~\ref{thm:stable_intervals_proper_hull}, the Helly hull $H(X)$ is locally finite. Since $G$ acts properly and cocompactly on $X$, we deduce that the action of $G$ on the Helly graph $H(X)$ is still proper and cocompact.
\ep

The notion of coarse injectivity can also be used to prove that a metric space is actually injective, as in the proofs of Theorems~\ref{thm:ell_infty_realization_poset_injective} and \ref{thm:local_global_injective_spaces}.

\bpro \label{pro:epsilon_coarse_injective}
Let $X$ denote a complete metric space that is $\eps$-coarsely injective for every $\eps>0$. Then $X$ is injective.
\epro

\bp We will prove that $X$ is injective: consider a family $(B_X(x_i,r_i))_{i \in I}$ of balls in $X$ such that $\forall i,j \in I, d(x_i,x_j) \leq r_i+r_j$. For any $\eps>0$, let us denote $A_\eps= \bigcap_{i \in I} B_X(x_i,r_i+\eps)$, which is non-empty by assumption of $\eps$-coarse injectivity.

Fix $0 < \eps \leq \eps'$, we will prove that the Hausdorff distance between $A_\eps$ and $A_{\eps'}$ is at most $\eps+\eps'$. Note that $A_\eps \subset A_{\eps'}$. Fix $x_0 \in A_{\eps'}$, we will prove that $d(x_0,A_\eps) \leq \eps+\eps'$. Assume that $0 \not\in I$, let $I_0=I \sqcup \{0\}$, and let $r_0=\eps'$. Consider the families $(x_i)_{i \in I_0}$ in $X$ and $(r_i)_{i \in I_0}$ in $\R_+$. For each $i,j \in I_0$, we know that $d(x_i,x_j) \leq r_i+r_j$: indeed, for any $i \in I$, we have $x_0 \in B_X(x_i,r_i+\eps')$. By $\eps$-coarse injectivity, we deduce that the intersection $\bigcap_{i \in I_0} B_X(x_i,r_i+\eps)$ is not empty. In particular, the ball $B_X(x_0,r_0+\eps)=B_X(x_0,\eps+\eps')$ intersects $A_\eps= \bigcap_{i \in I} B_X(x_i,r_i+\eps)$. This implies that $d(x_0,A_\eps) \leq \eps+\eps'$. So we have proved that the Hausdorff distance between $A_\eps$ and $A_{\eps'}$ is at most $\eps+\eps'$.

For each $n \in \N$, consider by induction $x_n \in A_{2^{-n}}$ such that for eacn $n \geq 0$ we have $d_X(x_{n+1},x_n) \leq 2^{-n}+2^{-(n+1)} \leq 2^{-n+1}$. For each $0 \leq n \leq m$, we have $d_X(x_n,x_m) \leq 2^{-n+2}$, hence the sequence $(x_n)_{n \in \N}$ is a Cauchy sequence in $X$. Since $X$ is complete, it has a limit $y \in X$. For each $n \in \N$, we have $y \in A_{2^{-n}}$, so for each $i \in I$ we have $d_X(y,x_i) \leq r_i+2^{-n}$. We deduce that, for each $i \in I$, we have $d_X(y,x_i) \leq r_i$. In other words, $y$ belongs to the intersection $\bigcap_{i \in I} B(x_i,r_i)$: we have proved that $X$ is injective.
\ep

\section{From lattices to injective metric spaces and Helly graphs} \label{sec:lattices_injective_helly}

\subsection{Posets and orthoscheme complexes} \label{subsec:orthoscheme}

\bdf[Poset, lattice, bowtie]
A \emph{poset} is a set endowed with some partial order.

A poset $(P,\leq)$ is called a \emph{meet-semilattice} if any $a,b \in P$ have a greatest lower bound $a \wedge b$, called the \emph{meet} of $a$ and $b$, i.e. such that
$$\forall c \in P, (c \leq a \mbox{ and } c \leq b) \Longrightarrow c \leq (a \wedge b).$$

Similarly, a poset $(P,\leq)$ is called a \emph{join-semilattice} if any $a,b \in P$ have a lowest upper bound $a \vee b$, called the \emph{join} of $a$ and $b$.

A poset $P$ is called a \emph{lattice} if it is both a meet-semilattice and a join-semilattice.

A \emph{bowtie} in a poset $(P,\leq)$ consists in four pairwise distinct elements $a,b \leq c,d$ such that $a$ and $b$ are not comparable, and neither are $c$ and $d$, and
$$\forall x \in P, a,b \leq x \leq c,d \Longrightarrow x \in \{a,b,c,d\}.$$

A \emph{chain} in a poset $P$ is a totally ordered subset $C \subset P$. If $C$ is finite, its \emph{length} is $|C|-1$.

A poset $(P,\leq)$ is \emph{homogeneous} if, for any $x \leq y$ in $P$, there is a bound on the lengths of chains from $x$ to $y$.

A poset $(P,\leq)$ is \emph{bounded} if it has a minimum element, denoted $0$, and a maximum element, denoted $1$.

If $(P,\leq)$ is a poset and $x \in P$, we will denote the subposets $P_{\leq x}=\{y \in X \st y \leq x\}$ and $P_{\geq x}=\{y \in X \st y \geq x\}$.
\edf

Bowties are the main obstructions to being a lattice, as is summarized in the following result.

\bpro
Let $P$ denote a homogeneous poset, and consider $P \cup \{0,1\}$ the poset with minimum $0$ and maximum $1$. The following are equivalent:
\bit
\item $P \cup \{0,1\}$ is a lattice.
\item $P \cup \{0\}$ is a meet-semilattice.
\item $P \cup \{1\}$ is a join-semilattice.
\item $P$ has no bowtie.
\eit
\epro

\bp
Assume that $L \cup \{0\}$ is a meet-semilattice, and consider pairwise distinct $a,b<c,d$ in $L$. Then the meet $x$ of $c,d$ is such that $a,b \leq x \leq c,d$. So either $x \not\in \{a,b,c,d\}$, or $a,b$ are comparable, or $c,d$ are comparable. Hence $L$ has no bowties.

\mk

Conversely, assume that $L$ has no bowtie. Note that $L \cup \{0\}$ has no bowtie either. Fix $a,b \in L$, and let $M$ denote the set of lower bounds of $a$ and $b$ in $L \cup \{0\}$: we have $0 \in M$, so $M$ is not empty. Let us consider a sequence $(x_n)_{n \in \N}$ in $M$ such that for each $n \in \N$, we have $x_n \leq x_{n+1}$. For each $n \in \N$, we have $x_0 \leq x_1 \leq \dots \leq x_n \leq a$, so it is a chain from $x_0$ to $a$. Since there is a bound on the length of such chains, we deduce that the sequence $(x_n)_{n \in \N}$ in $M$ is eventually constant.

\mk

We may therefore consider a maximal element $x$ of $M$, i.e. such that, for any $y \in M$, we have $x \not < y$. We will prove that $x$ is greatest in $M$, i.e. for any $y \in M$, we have $y \leq x$. It is sufficient to prove that $x$ is the unique maximal element in $M$: by contradiction, assume that $y \in M$ is a maximal element distinct from $x$. Then $x,y<a,b$ form a bowtie. Hence $x$ is the unique maximal element of $M$, and it is the meet of $a$ and $b$ in $L \cup \{0\}$.

Similarly, any two elements of $L$ have a meet in $L \cup \{1\}$. So $L \cup \{0,1\}$ is a lattice.
\ep

\bdf[Flag semilattice]
A meet-semilattice $(P,\leq)$ is called \emph{(upper) flag} if any $a,b,c \in P$ which are pairwise upperly bounded have an upper bound.

A join-semilattice $(P,\leq)$ is called \emph{(lower) flag} if any $a,b,c \in P$ which are pairwise lowerly bounded have a lower bound.
\edf

\brk
A bounded meet-semilattice is (upper) flag.
\erk

\bdf[Geometric realization]
Let $(P,\leq)$ denote a poset. Recall that the \emph{geometric realization} $|P|$ of $P$ is the simplicial complex with vertex set $P$, whose $k$-simplices are chains $v_0 < v_1 < \dots < v_k$ of length $k$. Note that simplices of $|P|$ have an induced total order on their vertices.
\edf

\subsection{Orthoscheme complexes} \label{subsec:orthoscheme_complexes}

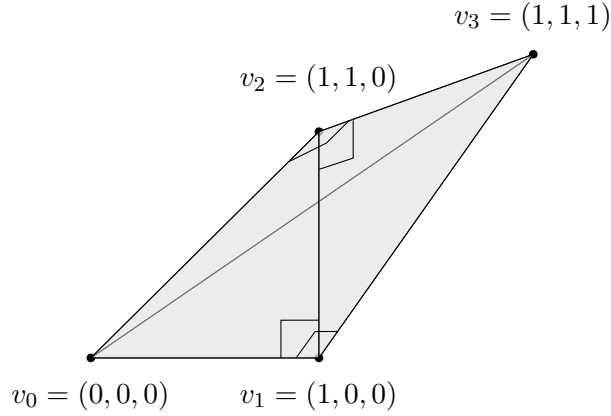
\begin{figure}
\begin{center}
\begin{tikzpicture}
\def \p {0.05}
\def \op {0.5}
\def \gris {black!10}
\draw[fill] (0,0) circle (\p) node(0) {};
\draw[fill] (3,0) circle (\p) node(1) {};
\draw[fill] (3,3) circle (\p) node(2) {};
\draw[fill] (3,3) + (20:3) circle (\p) node(3) {};

\draw[black,fill opacity=\op,fill=\gris] (0.center) -- (1.center) -- (3.center) -- (0.center);
\draw[black,fill opacity=\op,fill=\gris] (0.center) -- (2.center) -- (3.center) -- (0.center);
\draw[black,fill opacity=\op,fill=\gris] (0.center) -- (1.center) -- (2.center) -- (0.center);
\draw[black,fill opacity=\op,fill=\gris] (3.center) -- (1.center) -- (2.center) -- (3.center);
\node at ([yshift=-0.5cm]0) {\bfseries $v_0=(0,0,0)$};
\node at ([yshift=-0.5cm]1) {\bfseries $v_1=(1,0,0)$};
\node at ([yshift=0.7cm]2) {\bfseries $v_2=(1,1,0)$};
\node at ([yshift=0.5cm]3) {\bfseries $v_3=(1,1,1)$};

\draw (2.5,0) -- (2.5,0.5) -- (3,0.5);
\draw (2.7,0) -- (2.95,0.35) -- (3.25,0.35);
\draw (3,2.5) -- (3.45,2.65) -- (3.45,3.15);
\draw (2.6,2.6) -- (3.1,2.85) -- (3.4,3.15);

\end{tikzpicture}
\end{center}
\caption{The standard $3$-dimensional orthosimplex.}
\label{fig:orthosimplex}
\end{figure}

\bdf[Orthosimplex]
The \emph{standard orthosimplex} of dimension $n$ is the simplex of $\R^n$ with vertices $(0,\dots,0),(1,0,\dots,0),\dots,(1,1,\dots,1)$ (see Figure~\ref{fig:orthosimplex}). One may endow the simplex with the standard $\ell^p$ metric on $\R^n$, for any $p \in [1,\infty]$. For our purposes, we will only consider the $\ell^\infty$ metric. Note that any $n$-simplex with a total order on its vertices $v_0<v_1<\dots<v_n$ may be identified uniquely with the $\ell^\infty$ orthosimplex of dimension $n$, where each $v_i$ is identified with the vertex $(1,\dots,1,0,\dots,0)$ with $i$ ones and $n-i$ zeros. Also note that reversing the total order on the vertices gives rise to an isometry of the orthosimplex. 
\edf

\bdf[Orthoscheme realization]
Let $P$ denote a poset. Since each simplex of the geometric realization $X=|P|$ of $P$ has a total order on its vertices, one may endow each simplex with the standard $\ell^\infty$ orthosimplex metric. Let us then endow $X$ with the length metric associated with the piecewise $\ell^\infty$ metric on each simplex: it is called the \emph{$\ell^\infty$ orthoscheme metric} of $X$.
\edf

\bthm \label{thm:orthoscheme_complete_length}
Let $P$ denote a homogeneous poset. Then its $\ell^\infty$ orthoscheme realization $X$ is a complete length space.
\ethm

\bp
Note that since $X$ has finite dimension, the $\ell^\infty$ orthoscheme complex of $X$ has finitely many isometry types of cells: the standard orthoscheme $k$-simplices, where $k \leq n$. The proof of~\cite[Theorem~7.13]{bridson_haefliger} adapts without change to this situation.
\ep

We are now able to state the most interesting aspect of this construction. This result could be thought as the link condition for $\ell^\infty$ realizations of simplicial complexes.

\bthm[$\ell^\infty$ link condition] \cite[Theorem~6.3]{haettel_helly_kpi1} \label{thm:ell_infty_realization_poset_injective}
Let $P$ denote a homogeneous poset. Then its $\ell^\infty$ orthoscheme realization is locally injective if and only if, for every $x \in P$ :
\bit
\item The poset $P_{\geq x}$ is a flag meet-semilattice.
\item The poset $P_{\leq x}$ is a flag join-semilattice. 
\eit
Moreover, if $|P|$ is connected and simply connected, then $|P|$ is injective and admits a unique reversible, consistent, convex geodesic bicombing.
\ethm

\bp
The proof, quite technical, is presented in~\cite{haettel_helly_kpi1}. We just outline the main steps of the proof.

\mk

The first part is to prove the result when $P$ is a bounded homogeneous lattice. We will illustrate it when $P=\{0,x_1,x_2,x_3,1\}$, where $0 < x_1,x_2,x_3 < 1$ are the only relations. Then $|P|$ is the gluing of three two-dimensional orthosimplices along their longest edge $[0,1]$.

\mk

We first show that $|P|$ is a convex subset inside the gluing of the three copies $E_1,E_2,E_3$ of the half-plane $\{(x,y) \in \R^2 \st x \leq y\}$, when the $\ell^\infty$ metric: let us denote by $X=E_1 \cup E_2 \cup E_3$, where each half-plane is glued along the common boundary line $\{(x,y) \in \R^2 \st x=y\}$. It is then sufficient to prove that $X$ is injective.

\mk

We will use Theorem~\ref{thm:lattice_injective} to prove that $X$ is injective: there is a natural partial order on $X$ induced from the lattice order on each half-plane $E_1,E_2,E_3$, and there is a natural cofinal action of $\R$ on $X$ given by the diagonal action by translations. It remains to check that $X$ is a lattice.

\mk

For each $\eps>0$, we consider the graph $X_\eps$ which will be a discretization of $X$. It is the graph with vertex set $X_\eps=E_{1,\eps} \cup E_{2,\eps} \cup E_{3,\eps}$, where $E_{i,\eps}=E_i \cap (\eps \Z)^2$. There is an edge between two vertices of $X_\eps$ if their distance in $X$ is equal to $\eps$. We then prove that $X_\eps$, with the restriction of the partial order from $X$, is a lattice: this is the most technical part of the argument.

\mk

By considering the limit as $\eps$ goes to $0$, we are now able to prove that $X$ itself is a lattice. By applying Theorem~\ref{thm:lattice_injective}, we now deduce that $X$ is injective. As a consequence, the orthoscheme realization $|P|$ of the bounded lattice $P$ is injective.

\mk

In order to deal with the general case where $P$ is a flag semilattice, we use a discretization of $|P|$ as a graph, for which we prove the Helly property using the flag condition. We deduce that the orthoscheme realization $|P|$ is $\eps$-coarsely injective for each $\eps>0$, and we conclude by Proposition~\ref{pro:epsilon_coarse_injective} that $|P|$ is injective.
\ep

We can apply it to give another proof that CAT(0) cube complexes, endowed with the piecewise $\ell^\infty$ metric, are injective.

\bcor \label{cor:CCC_injective_semilattices}
Let $X$ denote a locally finite-dimensional CAT(0) cube complex. Then $(X,\ell^\infty)$ is injective.
\ecor

\bp
Let $P$ denote the poset consisting of all cubes of $X$, ordered by inclusion. Then the cubical barycentric subdivision of $X$ naturally identifies with the geometric realization of $P$: in particular, it is connected and simply connected.

\mk

Since $X$ is locally finite-dimensional, $P$ is homogeneous. Let $Q \in P$ denote an arbitrary cube in $X$.

The poset $P_{\leq Q}$ is the Boolean lattice $\{0,1\}^{\dim Q}$, which is in particular a flag join-semilattice.

The link of $Q$ in $X$ is simplicial if and only if the poset $P_{\geq Q}$ is a meet-semilattice.

The link of $Q$ in $X$ is a flag complex if and only if that the poset $P_{\geq Q}$ is flag.
\ep

\bexe (Exercise)
Let $X$ denote a Helly graph, and let
$$P_X = \{\mbox{ round cliques of $X$}\} \subset {\cal P}(X),$$
with the inclusion partial order. Show that $P_X$ satisfies the assumptions of Theorem~\ref{thm:ell_infty_realization_poset_injective}.
\eexe

\subsection{Lattices, injective metrics and Helly graphs}

If $L$ is a poset, we say that a non-empty subset $A \subset L$ has a join if $L_{\geq A}=\{x \in L \st \forall a \in A, a \leq x\}$ has a lowest element.

\bthm \label{thm:lattice_injective}
Let $L$ denote a lattice such that every upperly bounded subset of $L$ has a join. Let $H=\Z$ or $\R$, and assume that we have an order-preserving, continuous (with respect to the order topology on $L$), cofinal action of $H$ on $L$ (denoted $+$), i.e.
$$\forall x,y \in L, \exists t \in H \cap \R_+, x-t \leq y \leq x+t.$$
For each $x,y \in L$, let us define
$$d(x,y)=\inf\{t \in H \cap \R_+, x-t \leq y \leq x+t\}.$$
Then we have the following:
\bit
\item If $H=\R$, then $(L,d)$ is an injective metric space.
\item If $H=\Z$, then $(L,d)$ is the vertex set of a Helly graph, with the combinatorial metric.
\eit
\ethm

\bexe
Let us consider $L=\R^n$, with the following order: $x \leq y$ if and only if $\forall 1 \leq i \leq n, x_i \leq y_i$. It is clear that $L$ is a lattice, and that every upperly bounded subset has a join. Let us consider the action of $H=\R$ on $\R^n$ by $t \cdot x=(x_1+t,x_2+t,\dots,x_n+t)$: this action is cofinal. The metric from Theorem~\ref{thm:lattice_injective} is the standard $\ell^\infty$ metric on $\R^n$.

If we restrict to $L=\Z^n$ and $H=\Z$, we obtain the standard Helly graph with vertex set $\Z^n$.
\eexe

\bexe (Exercise)
The gluing of two copies of the plane $(\R^2,\ell^\infty)$ along the half-plane $\{(x,y) \in \R^2 \st y \leq x\}$ is injective.
\eexe

\bp
It is not hard to see that $d$ is symmetric and satisfies the triangle inequality. The continuity assumption of the action of $H$ ensures that $d$ is positive.

\mk

Note that, by definition, for any $x \in L$ and any $r \in H_+$, we have the equality
$$B(x,r) = \{y \in X \st x-r \leq y \leq x+r\} = I(x-r,x+r)$$
so balls in $L$ are intervals.

Moreover, by assumption on $L$, it is easy to see that intervals in $L$ satisfy the Helly property: fix any collection $(I(x_i,y_i))_{i \in I}$ of pairwise intersecting intervals of $L$. Fix $j_0 \in I$. For any $i \in I$, we have $x_i \leq y_{j_0}$, so the set $\{x_i, i \in I\}$ is upperly bounded. By assumption, the set $\{x_i, i \in I\}$ has a join $z \in L$, such that $z \leq y_{j_0}$. This holds for any $j_0 \in I$, so $z$ belongs to the intersection of all intervals $(I(x_i,y_i))_{i \in I}$.

\mk

According to Theorem~\ref{thm:injective_equivalent_pairwise_intersecting_balls}, it is enough to prove that $(L,d)$ is a geodesic metric space (if $H=\R$) or that $(L,d)$ is the vertex set of a connected graph (if $H=\Z$).

\mk

Consider the case $H=\Z$. Let us define the graph $\Gamma$ with vertex set $L$, and with an edge between $x, y \in L$ if $x-1 \leq y \leq x+1$: we will prove that the graph $\Gamma$ is connected. We will first prove, by induction on $k \geq 0$, that for any $x \in L$, the ball $B_{\Gamma}(x,k)$ in the graph $\Gamma$ coincides with the interval $I(x-k,x+k)=\{y \in L \st x-k \leq y \leq x+k\}$. 

\mk

For $k \leq 1$ it is the definition of the edges of $\Gamma$, so fix $k \geq 2$ and assume that the statement holds for $k-1$. Fix $y \in I(x-k,x+k)$, we will prove that $y \in B_{\Gamma}(x,r)$. Since $y \geq x-k$, we deduce that $y+1 \geq x-k+1$, and also since $y \leq x+k$ we deduce that $y-1 \leq x+k-1$. So we have that both $y+1$ and $x+k-1$ are superior to both $y-1$ and $x-k+1$: since $L$ is a lattice, there exists some element $z \in L$ in the intersection $I(y-1,y+1) \cap I(x-k+1,x+k-1)$. In particular, $y$ and $z$ are adjacent in $\Gamma$, and by induction we know that $d_{\Gamma}(z,x) \leq k-1$, so $d_{\Gamma}(x,y) \leq k$. Conversely, it is clear that the ball $B_{\Gamma}(x,k)$ is included in the interval $I(x-k,x+k)$. So we have $B_{\Gamma}(x,k) = I(x-k,x+k)$.

\mk

Since the action of $H$ is cofinal on $L$, we deduce that the graph $\Gamma$ is connected, and furthermore that $d_{\Gamma}=d$. Since balls in $(L,d)$ satisfy the Helly property, we conclude that $(L,d)$ is the vertex set of the Helly graph $\Gamma$.

\mk

Consider the case $H=\R$: fix $\eps>0$, we will first prove that $(L,d)$ is $\eps$-coarsely injective. Let us consider the action of $H_\eps=\eps\Z \subset \R$ on $L$, and the associated Helly graph distance
$$\forall x,y \in L, d_\eps(x,y) = \inf\left\{t \in \eps \N \st x-t \leq y \leq x+t\right\}.$$
It is easy to see that $d \leq d_\eps \leq d+\eps$. Since ball for $d_\eps$ satisfy the Helly property, we deduce that $(L,d)$ is $\eps$-coarsely injective, and this holds for all $\eps>0$.

We will now prove that $(L,d)$ is complete: let $(x_n)_{n \in \N}$ denote a Cauchy sequence in $L$. Since $(d(x_0,x_n))_{n \in \N}$ is bounded, the sequence $(x_n)_{n \in \N}$ has an upper bound and a lower bound. For each $n \in \N$, let us denote
$$y_n=\bigvee_{k \geq n} x_k.$$
For each $0 \leq n \leq k$, we have $d(y_n,x_k) \leq \sup_{k,k' \geq n} d(x_k,x_{k'})$. The sequence $(y_n)_{n \in \N}$ is decreasing and is bounded below: let us define
$$z=\bigwedge_{n \in \N} y_n.$$
We deduce that for each $n \in \N$, we have
$$d(y_n,z) \leq \sup_{k \geq n} d(y_n,x_k) \leq \sup_{k,k' \geq n} d(x_k,x_{k'}).$$
Since $(x_n)_{n \in \N}$ is a Cauchy sequence, we conclude that $(y_n)_{n \in \N}$ converges to $z$, and hence $(x_n)_{n \in \N}$ converges to $z$. So $(L,d)$ is complete.

According to Proposition~\ref{pro:epsilon_coarse_injective}, we conclude that $(L,d)$ is injective.
\ep

\subsection{Application to symmetric spaces}

We will explain how we can use Theorem~\ref{thm:lattice_injective} to study the injective hull of the symmetric space of $\GL(n,\R)$. More details can be found in~\cite{haettel_injective_buildings}.

\mk

Fix $n \geq 2$. Let us recall that the symmetric space $X_0$ of $\SL(n,R)$ is the manifold
$$X_0=\SL(n,\R)/\SO(n),$$
endowed with the unique (up to homothety) $SL(n,\R)$-equivariant Riemannian metric: it is of nonpositive sectional curvature (see for instance~\cite{helgason}, \cite{eberlein}, \cite{guivarch} and \cite{bridson_haefliger}). Since $\SO(n)$ is the stabilizer in $\SL(n,\R)$ of the standard ellipsoid, or equivalently of the standard Euclidean norm on $\R^n$, we have the following models for $X_0$:
\beq X_0& \simeq &\{\mbox{volume  $1$ ellipsoids in $\R^n$}\}\\ 
&\simeq& \{\mbox{Euclidean norms on $\R^n$}\} / \mbox{homothety},\eeq
where by \emph{ellipsoid} we mean a centered (at $0$), non-degenerate ellipsoid. Indeed, one may associate to each Euclidean norm on $\R^n$ its unit ball, which is a volume $1$ ellipsoid in $\R^n$ (up to homothety).

\mk

The symmetric space $X$ of $\GL(n,\R)$ is the manifold
$$X=\GL(n,\R) / O(n) \simeq X_0 \times \R,$$
it is also a Riemannian manifold with nonpositive sectional curvature, with a flat factor $\R$. The space $X$ can equivalently be described as:\beq X& \simeq &\{\mbox{ellipsoids in $\R^n$}\}\\ 
&\simeq& \{\mbox{Euclidean norms on $\R^n$}\}.\eeq
The space $X$ is naturally a poset for any of these two models: say that two ellipsoids $C,C' \subset \R^n$ satisfy $C \leq C'$ if and only if $C \subset C'$. Equivalently, say that two Euclidean norms $N,N'$ on $\R^n$ satisfy $N \leq N'$ if and only if $\forall v \in \R^n, N(v) \geq N'(v)$. Indeed, having a smaller unit ball is equivalent to having a larger norm.

\mk

If $n \geq 2$, the poset $(X,\leq)$ is not a lattice: indeed, consider two generic ellipses $C,C'$ in $\R^2$. There are many maximal ellipses contained in the intersection $C \cap C'$. This remark motivates the following introduction: let us denote by $\hat{X}$ the space
\beq \hat{X}& \simeq &\{\mbox{symmetric, closed, with non-empty interior, convex bodies in $\R^n$}\}\\ 
&\simeq& \{\mbox{norms on $\R^n$}\}.\eeq
The space $\hat{X}$ also has a natural partial order. However, we now have the following.

\bpro
The poset $\hat{X}$ is a lattice. Moreover, any upperly bounded subset $A \subset \hat{X}$ has a join. 
\epro

\bp
Let $C,C' \in \hat{X}$ any convex bodies in $\R^n$. One easily checks that $C \cap C'$ is the meet of $C$ and $C'$, and that the convex hull $\Hull(C \cup C')$ is the join of $C$ and $C'$. Hence $\hat{X}$ is a lattice.

\mk

Now let $A \subset \hat{X}$ denote a non-empty upperly bounded subset: there exists $C_0 \in \hat{X}$ such that $\forall C \in A, C \subset C_0$. Let us now consider $D \in \hat{X}$ the convex hull of $\cup A$, i.e. the smallest closed convex subset of $C_0$ containing every element of $A$. One sees that $D$ is the join of $A$.  
\ep

There is a natural action of $\R$ on $\hat{X}$ (which extends the translation action on $X=X_0 \times \R$):
$$\forall t \in \R, \forall C \in \hat{X}, t \cdot X = e^t C.$$
It is easy to see that this action is cofinal, so we deduce that Theorem~\ref{thm:lattice_injective} provides a metric $d$ on $\hat{X}$ such that $(\hat{X},d)$ is injective. Note that $\GL(n,R)$ has a natural action on $\hat{X}$, by the standard linear action on convex bodies on $\R^n$ or on norms on $\R^n$.

\blem
The metric $d$ on $\hat{X}$ is invariant under the action of $\GL(n,\R)$.
\elem

\bp
Let $C,C'$ denote convex bodies in $\R^n$, and assume that $t>0$ is such that $e^{-t}C \subset C' \subset e^t C$. Then, for any $g \in \GL(n,\R)$, one has
$$g \cdot (e^{-t}C)=e^{-t} (g \cdot C) \subset g \cdot C' \subset g \cdot (e^{t}C)=e^{t} (g \cdot C).$$
One easily deduces that $d(g \cdot C,g \cdot C')=d(C,C')$.
\ep

\brk
When restricted to the symmetric space $X$, the metric $d$ is quite simple to describe: any two points in $X$ lie in a common flat, so it is sufficient to describe the metric $d$ on the standard flat $A \cdot O(n)$, where $A \subset \GL(n,R)$ is the diagonal subgroup with entries in $(0,\infty)$. Then the restriction of the metric $d$ to $A \simeq \R^n$ is the standard $\ell^\infty$ metric on $\R^n$.

In fact, one can prove that $(\hat{X},d)$ is precisely the injective hull of $(X,d)$ (see~\cite[Theorem~F]{haettel_injective_buildings}).
\erk

In fact, one can say more about the action of $\GL(n,\R)$ on $\hat{X}$.

\bpro
The metric space $\hat{X}$ is proper, and the isometric action of $\GL(n,\R)$ is proper and cocompact.
\epro

\bp
Let $K \in \hat{X}$. Let $B \subset K$ denote the unique John-L\"owner ellipsoid of maximal volume: we have $B \in X$. According to~\cite{john}, we know that $\hat{d}(B,K) \leq \log(\sqrt{n})$. Therefore any point of $\hat{X}$ is at distance at most $\log(\sqrt{n})$ from $X$, on which $\GL(n,\R)$ acts transitively. So the isometric action of $\GL(n,\R)$ on $\hat{X}$ is cobounded.

Given a convex body $C \in \hat{X}$ and $t > 0$, it is easy to see that the set $B(C,t)$ of all convex bodies $C' \in \hat{X}$ such that $e^{-t}C \subset C' \subset e^tC$ is compact for the Hausdorff topology. Since the metric $d$ defines the same topology, we deduce that the ball $B(C,t)$ is compact. Hence $\hat{X}$ is a proper metric space.

Since the action of $\GL(n,\R)$ on $X$ is proper, we conclude that the action of $\GL(n,\R)$ on $\hat{X}$ is proper and cocompact.
\ep

\bcor
Any uniform lattice in $\GL(n,\R)$ has a proper and cocompact action by isometries on an injective metric space.
\ecor

\brk
Non-uniform lattices in $\GL(n,\R)$, with $n \geq 3$, have distorted elements so they cannot act properly and cocompactly (even coboundedly) by isometries on an injective metric space.
\erk

\bq
Are uniform lattices in $\GL(n,\R)$ Helly? Note that biautomaticity is open for such groups.
\eq

\bq
Which uniform lattices in semisimple Lie groups are injective? There are partial answers in~\cite{haettel_injective_buildings}.
\eq

Many other classical symmetric spaces may be realized as fixed-point set of involutions inside the symmetric space of $\GL(n,\R)$, so one can also deduce their isometry groups act properly and cocompactly on an injective metric space: see~\cite[Theorem~G]{haettel_injective_buildings}.

\bq
What is the injective hull $E(\H^2)$ of the hyperbolic plane $\H^2=\SL(2,\R)/\SO(2)$? Note that $E(\H^2)$ can naturally be described of the injective metric space $\hat{X}$ consisting of all convex bodies or $\R^2$, and according to Proposition~\ref{pro:hull_hyperbolic_space} it is at bounded distance from the subspace $\H^2$ consisting of all volume $1$ ellipses.
\eq

\subsection{Application to buildings}

We will explain how we can use Theorem~\ref{thm:lattice_injective} to study the Helly hull of the Bruhat-Tits building of $\GL(n,\K)$, where $\K$ is a non-Archimedean local field. More details can be found in~\cite{haettel_injective_buildings}.

\mk

Fix $n \geq 2$. Fix a non-Archimedean local field $\K$, i.e. $\K$ is a finite extension of either the field $\Q_p$ of $p$-adic numbers, for some prime number $p$, or a finite extension of the field $\F_q((t))$ of Laurent series, where $\F_q$ is a finite field. Let $q$ denote the cardinality of the residue field, and consider the absolute value $|\cdot| : \K \ra q^\Z \cup \{0\} \subset \R_+$. Let us recall that a map $\eta : \K^n \ra \R_+$ is an \emph{ultrametric norm} on $\K^n$ if it satisfies the following.
\bit
\item $\forall v \in V, \eta(v)=0 \Longleftrightarrow v=0$.
\item $\forall v \in V, \forall \alpha \in \K, \eta(\alpha v) = |\alpha| \eta(v)$.
\item $\forall u,v \in V, \eta(u+v) \leq \max(\eta(u),\eta(v))$.
\eit

\mk

Let us recall that the Bruhat-Tits building of $\SL(n,\K)$ has a natural simplicial structure, whose vertex set is
$$X_0=\{\mbox{ultrametric norms $\K^n \ra q^\Z \cup \{0\}$}\} / \mbox{homothety by $q^\Z$}.$$
Similarly to the case of symmetric spaces, we will consider the Bruhat-Tits building of $\GL(n,\K)$: it has a natural simplicial structure, whose vertex set is
$$X=\{\mbox{ultrametric norms on $\K^n \ra q^\Z \cup \{0\}$}\}.$$
This space has been defined and studied by Goldman and Iwahori before the invention of buildings, see~\cite{goldman_iwahori}.

Note that $X$ is a poset, where we say that two ultrametric norms $N,N'$ on $\K^n$ satisfy $N \leq N'$ if and only if $\forall v \in \K^n, N(v) \geq N'(v)$. Contrary to the case of symmetric spaces, the space $X$ is already a lattice, so the situation is simpler:

\bpro
The poset $X$ is a lattice, and each upperly bounded subset has a join.
\epro

There is a natural action of $H=\Z$ on $X$ by homothety:
$$\forall t \in \Z, \forall N \in X, t \cdot N = q^t C.$$

It is easy to see that this action is cofinal, so we deduce that Theorem~\ref{thm:lattice_injective} provides a metric $d$ on $X$ whose underlying graph $\Gamma$ is Helly. Note that $\GL(n,\K)$ has a natural action on $X$, by the standard linear action by precomposition on ultrametric norms on $\K^n$. This action preserves the metric $d$, so we deduce the following.

\bpro
The group $\GL(n,\K)$ acts properly and cocompactly by automorphisms on the locally finite Helly graph $\Gamma$, which has the same vertex set as the Bruhat-Tits building of $\GL(n,\K)$.
\epro

\bcor
Any uniform lattice in $\GL(n,\K)$ is Helly, i.e. it has a proper and cocompact action by automorphisms on a Helly graph.
\ecor

Many other classical Euclidean buildings may be realized as fixed-point set of finite groups inside the Bruhat-Tits building of $\GL(n,\K)$, so one can also deduce their isometry groups act properly and cocompactly on a Helly graph: see~\cite[Corollary~D]{haettel_injective_buildings}.

\subsection{Application to Garside groups and Artin groups} \label{subsec:garside}

The notion of Garside groups originated in Garside's work on word and conjugacy problems for braid groups (\cite{garside_braid_groups}). The central idea is that braid groups can be endowed with a specific partial order which is a lattice. Garside groups then have been defined formally and studied by Dehornoy and Paris (see~\cite{dehornoy_paris_gaussian}), and in more depth in~\cite{garside}.

\mk

Here is one definition of Garside groups. We refer the reader to~\cite{garside} and \cite{mccammond_intro_garside} for more background on Garside groups. 

\bdf[Garside group]
\label{def:garside}
Let $G$ denote a group, $S \subset G$ a finite subset and $\Delta \in G$. The triple $(G,S,\Delta)$ is called a \emph{Garside structure} if the following conditions hold. Let $G^+$ denote the submonoid of $G$ generated by $S$.
\ben
\item The group $G$ is generated by $S$.
\item For any element $g \in G^+$, there is a bound on the length $n$ of expressions $g=s_1 \dots s_n$, where $s_1,\dots,s_n \in S \bs \{1\}$.
\item We define the partial $\leq_L$, $\leq_R$ on $G^+$ by $a\leq_L b$ if and only if $b=ac$ for some $c\in G^+$ and $a\leq_R b$ if and only if $b=ca$ for some $c\in G^+$.
The left $\leq_L$ and right $\leq_R$ orders on $G^+$ are lattices.
\item The set $S$ is a balanced interval between $1$ and $\Delta$, i.e.
$$S = \{g \in G^+ \st 1 \leq_L g \leq_L \Delta\} =  \{g \in G^+ \st 1 \leq_R g \leq_R \Delta\}.$$
\een
A group is called \emph{Garside} if it admits such a Garside structure, and $\Delta$ is called the \emph{Garside element}.
If the set $S$ is allowed to be infinite, we may say that $(G,S,\Delta)$ is a \emph{quasi-Garside structure}.
\edf

In~\cite{haettel_huang_weakly_modular}, we give equivalent characterizations of Garside groups, which are more geometric in flavour. We need first to define what Garside lattices and Garside flag complexes are.

\bdf[Garside lattice]
A \emph{Garside lattice} is a pair $(L,\varphi)$, where $L$ is a homogeneous lattice and $\varphi$ is an increasing automorphism of $L$, such that, for any $x,y \in L$, there exists $k \in \N$ such that $x \leq \varphi^k(y)$.
\edf

\bdf[Garside flag complex]
A \emph{Garside flag complex} is a pair $(X,\varphi)$, where $X$ is a simply connected flag simplicial complex with finite simplices, with a consistent total order on each simplex, and $\varphi$ is an order-preserving automorphism of $X$, such that the following hold:
\bit
\item For any simplex $\sigma$ of $X$, we have that $\sigma \cup \varphi(\min \sigma)$ is a simplex of $X$.
\item For any vertex $x \in X$, we have $\varphi(x)>x$, and the interval $[x,\varphi(x)]$ is a homogeneous lattice.
\eit
\edf

We are now able to state alternative characterizations of Garside groups.

\bthm \cite[Theorem~4.7]{haettel_huang_weakly_modular} \label{thm:dictionary_garside_group}
Let $G$ denote a group. The following are equivalent:
\ben
\item $G$ is a Garside group.
\item There exists a Garside lattice $(L,\varphi)$ such that $G$ can be realized as a group of order-preserving automorphisms of $L$ commuting with $\varphi$, acting freely and transitively on elements of $X$.
\item There exists a Garside flag complex $(X,\varphi)$ such that $G$ can be realized as a group of order-preserving automorphisms of $X$ commuting with $\varphi$, acting freely and transitively on vertices of $X$.
\een
\ethm

\bp
We will only briefly describe how one passes from one structure to the other.
\bit
\item[$1. \Rightarrow 2.$] Let $(G,S,\Delta)$ denote a Garside structure on $G$. Let $L=G$, with order $\leq_L$ defined by $a \leq_L b$ if and only if there exists $h \in G^+$ such that $b=ah$. Then $\varphi$ is the automorphism of $L$ given by the right multiplication by $\Delta$, and $G$ acts on $L$ by left multiplication.
\item[$2. \Rightarrow 3.$] Let $X$ denote the flag simplicial complex with vertex set $L$, with a $k$-simplex for each chain $x_0 < x_1 < \dots < x_k$ such that $x_k \leq \varphi(x_0)$.
\item[$3. \Rightarrow 1.$] Since $G$ acts freely and transitively on vertices of $X$, we may equivariantly identify the vertices of $X$ with $G$. Now $\Delta$ corresponds to the element $\varphi(e)$, and $S$ consists in all elements $g \in G$ such that $\{e,g,\Delta\}$ is a simplex of $X$.
\eit
\ep

The main example of Garside groups are braid groups, and more generally spherical type Artin groups, which we briefly define now. We refer the reader to~\cite{kassel_turaev} and \cite{farb_margalit} for references concerning braid groups, and to~\cite{paris_kpi1,godelle_paris} for references concerning Artin groups.

\bdf[Coxeter and Artin group]
Consider a finite simple graph $\Gamma$ with vertex set $S$ and with edges labeled by some integer in $\{2,3,\dots\}$. One associates the \emph{Coxeter group} $W(\Gamma)$ with the following presentation:
\beq W(\Gamma) = \<S &\st & \forall s \in S, s^2=1,\\
&&\forall \{s,t\} \in \Gamma^{(1)}, [s,t]_m=[t,s]_m \mbox{ if the edge $\{s,t\}$ is labeled $m$}\>,\eeq
where $[s,t]_m$ denotes the word $ststs\dots$ of length $m$.

One also associates the \emph{Artin group} $A(\Gamma)$ with the following presentation:
$$A(\Gamma) = \<S \st \forall \{s,t\} \in \Gamma^{(1)}, [s,t]_m=[t,s]_m \mbox{ if the edge $\{s,t\}$ is labeled $m$}\>.$$

If the Coxeter group $W(\Gamma)$ is finite, the Artin group $A(\Gamma)$ is called of \emph{spherical type}.
\edf

\bexe\

\bit
\item If $\Gamma$ is the graph with no edge, then $W(\Gamma) = \Z/2\Z \star \dots \star \Z/2\Z$ and $A(\Gamma)=\F(S)$.
\item If $\Gamma$ is the complete graph with all edges labeled $2$, then $W(\Gamma) = (\Z/2\Z)^S$ and $A(\Gamma)=\Z^S$.
\item Let $\Gamma$ denote the complete graph on $S=\{s_1,\dots,s_{n-1}\}$, with all edges labeled $2$ except the edges between $s_i$ and $s_{i+1}$ labeled $3$, for $1 \leq i \leq n-2$. Then $W(\Gamma)$ is the symmetric group $\frak{S}_n$ and $A(\Gamma)$ is the $n$-strand braid group. 
\eit
\eexe

The following is key in the study of spherical type Artin groups.

\bthm \label{thm:spherical_artin_garside} \cite{deligne,brieskorn_saito,bessis,birman_ko_lee}
Any spherical type Artin group is a Garside group.
\ethm

We are now able to state the main result relating Garside groups, injective metric spaces and Helly graphs.

\bthm \cite{huang_osajda_helly,haettel_helly_kpi1} \label{thm:garside_implies_helly}\

Any Garside group is a Helly group.

More precisely, if $(G,S,\Delta)$ is a Garside structure, then the Cayley graph of $G$ with respect to $S^{-1}S$ is Helly.

Moreover, if $X$ denotes a Garside flag complex on which $G$ acts freely transitively as in Theorem~\ref{thm:dictionary_garside_group}, then the $\ell^\infty$ orthoscheme realization of $X$ is injective.
\ethm

\bp
The first proof, due to Huang and Osajda (see~\cite{huang_osajda_helly}), uses the notion of cell Helly complexes. Let $(G,S,\Delta)$ denote a Garside structure on $G$. Let $X$ denote the cell complex with vertex set $S$, and with cells $gS$, for all $g \in G$. Then we admit that the complex $X$ is cell Helly. According to Theorem~\ref{thm:thickening_helly}, this implies that the thickening of $X$, which coincides with the Cayley graph $\Gamma$ of $G$ with respect to $S^{-1}S$, is Helly.

\mk

Another proof (see~\cite{haettel_helly_kpi1}) uses cofinal actions on lattices. For instance, let us 
a Garside lattice $(L,\varphi)$ on which $G$ acts freely transitively as in Theorem~\ref{thm:dictionary_garside_group}. Then the action of $H=\Z$ generated by $\varphi$ is cofinal on $L$. We can then apply Theorem~\ref{thm:lattice_injective} and deduce that the corresponding graph, which turns out to be the Cayley graph $\Gamma$, is Helly.

\mk

The last statement concerning the orthoscheme realization of the Garside flag complex is an immediate application of Theorem~\ref{thm:ell_infty_realization_poset_injective}.
\ep

\brk
In particular, a consequence of Theorem~\ref{thm:garside_implies_helly} is that the dual braid complex $X_n$ for the $n$-strand braid group, endowed with the $\ell^\infty$ orthoscheme metric, is injective. However, the question whether $X_n$, endowed with the $\ell^2$ orthoscheme metric, could be CAT(0) is still open (see~\cite{brady_mccammond,b6}). In particular, the question whether braid groups are CAT(0) is still open, but we do know that they are Helly.
\erk

\bq
Is every Helly group CAT(0)?
\eq

\section{Summary of examples and properties} \label{sec:summary_examples_properties}

In this section, we give a summary of most known Helly and injective groups, as well as the main consequences.

\bthm The following groups are Helly, i.e. admit proper and cocompact ations by isometries on Helly graphs:
\bit
\item Hyperbolic groups (\cite{lang}).
\item Relatively hyperbolic groups with Helly parabolics (\cite{osajda_valiunas}).
\item Cocompactly cubulated groups.
\item Braid groups, and more generally FC type Artin groups (\cite{huang_osajda_helly}).
\item Garside groups (\cite{huang_osajda_helly}).
\item Uniform lattices in $\GL(n,\K)$, where $\K$ is a non-Archimedean local field, or in other classical groups (\cite{haettel_injective_buildings}).
\item Few crystallographic groups (\cite{hoda:crystallographic}).
\eit

The following groups are injective, i.e. admit proper and cocompact ations by isometries on injective metric spaces:
\bit
\item Uniform lattices in $\GL(n,\R)$, or in other classical groups (\cite{haettel_injective_buildings}).
\item Uniform lattices in products of Gromov-hyperbolic spaces (and there exists a non-Helly example in the product of $\H^2$ and a tree, see~\cite{hugues_valiunas_HHS}).
\eit

The following groups are coarsely injective, i.e. admit proper and cobounded ations by isometries on injective metric spaces:
\bit
\item Mapping class groups of finite type surfaces (\cite{haettel_hoda_petyt}).
\item More generally, hierarchically hyperbolic groups (\cite{haettel_hoda_petyt}).
\eit
\ethm

Many stability properties are true for Helly groups, see~\cite[Theorem~1.3]{helly_groups}, notably the stability by graph products. Note that similar results probably hold for injective groups and coarsely injective groups.

\mk

Let us now gather consequences for a group to be Helly, injective, or coarsely injective. Recall that any Helly group is an injective group, and every injective group is coarsely injective.

\bthm\

\noindent If $G$ is a coarsely injective group, the following hold:
\bit
\item The group $G$ has finitely many conjugacy classes of finite subgroups (\cite[Proposition~1.2]{lang}).
\item The group $G$ has at most Euclidean Dehn functions (\cite{wenger}).
\item The group $G$ has contractible asymptotic cones (\cite[Theorem~1.5]{helly_groups}).
\item The group $G$ satisfies the coarse Baum-Connes conjecture (\cite{fukaya_oguni}).
\item The group $G$ is strongly shortcut (\cite[Theorem~B]{haettel_hoda_petyt}).
\item The group $G$ is semihyperbolic (see~\cite{lang}), and in particular (\cite{bridson_haefliger}):
\bit
\item The centralizer of a finite set of elements of $G$ is semihyperbolic.
\item $G$ has solvable word and conjugacy problems.
\item Any polycyclic subgroup of $G$ is virtually abelian, finitely generated and undistorted.
\eit
\eit
If $G$ is an injective group, the following extra consequences hold:
\bit
\item The group $G$ satisfies the Farrell-Jones conjecture with finite wreath products (\cite[Theorem~6.1]{kasprowski_rueping}).
\item The group ring $\K[G]$ satisfies Kaplansky's idempotent conjecture, if $\K$ is a field with characteristic zero (\cite[Theorem~0.12]{bartels_luck_holger_farrell_jones}).
\eit
If $G$ is a Helly group, the following extra consequences hold:
\bit
\item $G$ acts properly and cocompactly on a contractible finite-dimensional simplicial complex (\cite{lang}).
\item Any element of $G$ has rational translation length, with uniformly bounded denominator (\cite[Theorem~O]{haettel_osajda_locally_elliptic}).
\item The group $G$ is biautomatic (\cite[Theorem~1.5]{helly_groups}).
\item Any torsion subgroup of $G$ is finite (\cite[Corollary~I]{haettel_osajda_locally_elliptic}).
\eit
\ethm

\section{Overview of actions on $L^1$, $L^2$ and $L^\infty$ spaces} \label{sec:L1L2Linfini}

Let us now give a very brief and informal overview of examples of groups acting "nicely" on "$L^p$ spaces", with some precisions:
\bit
\item A "$L^1$ space" means a median graph (i.e. a CAT(0) cube complex) or a metric median space.
\item A "$L^2$ space" means a CAT(0) space.
\item A "$L^\infty$ space" means a Helly graph or an injective metric space.
\item A "nice action on a $L^p$ space" means a proper action by isometries, which is assumed to be cocompact or cobounded  for $L^2$ or $L^\infty$ spaces.
\eit
We insist that the following tables are to be considered as coarse summaries. 

\mk

\begin{tabular}{|l|l|l|l|}
\hline
Groups acting nicely on... & $L^1$ & $L^2$ & $L^\infty$ \\
\hline
Hyperbolic groups & no & open & yes\\
\hline
Cubulable groups & yes & yes & yes\\
\hline
Relatively hyperbolic groups & no & open & yes\\
\hline
Coxeter groups & yes & yes & no\\
\hline
Braid groups & open & open & yes\\
\hline
Garside groups & no & open & yes\\
\hline
Artin groups & no & open & open\\
\hline
Mapping class groups of surfaces & open & no & yes\\
\hline
Higher rank lattices & no & yes & many\\
\hline
\end{tabular}

\mk

Let us give a similar overview of properties of groups acting nicely on such spaces.

\mk

\begin{tabular}{|l|l|l|l|}
\hline
Properties of a group acting nicely on... & $L^1$ & $L^2$ & $L^\infty$ \\
\hline
Finite number of conjugacy classes of finite subgroups & yes & yes & yes\\
\hline
Semihyperbolicity & yes & yes & yes\\
\hline
Biautomaticity & yes & no & yes\\
\hline
Euclidean isoperimetric inequalities & yes & yes & yes\\
\hline
No infinite torsion subgroup & yes & open & yes\\
\hline
Tits alternative & yes & open & open\\
\hline
Rank rigidity & yes & open & open\\
\hline
Strong bolicity & yes & yes & no\\
\hline
Property RD & yes & open & open\\
\hline
Farrell-Jones conjecture & yes & yes & yes\\
\hline
Baum-Connes conjecture & yes & open & open\\
\hline
\end{tabular}

\section{Open questions} \label{sec:questions}

We end with a list of open questions emphasizing the fact that these topics are very active.

\ben
\item What is the Helly hull of the standard Cayley graph of a Coxeter group?
\item Which Coxeter groups are Helly?
\item Which Artin groups are Helly? As a simple example, is the $\tilde{A}_2$ Artin group Helly? According to~\cite{haettel_huang_garside_artin_product_Z}, its direct product with $\Z$ is Helly.
\item What is the injective hull $E(\H^2)$ of the hyperbolic plane $\H^2=\SL(2,\R)/\SO(2)$? Note that $E(\H^2)$ can be described inside the injective metric space $\hat{X}$ consisting of all convex bodies or $\R^2$, and according to Proposition~\ref{pro:hull_hyperbolic_space} it is at bounded distance from the subspace $\H^2$ consisting of all volume $1$ ellipses.
\item What is the injective hull of the hyperbolic $n$-space?
\item Are uniform lattices in $\GL(n,\R)$ Helly? Note that biautomaticity is open for such groups.
\item Which uniform lattices in semisimple Lie groups are injective? There are partial answers in~\cite{haettel_injective_buildings}.
\item Do Helly groups satisfy a Tits alternative?
\item Are Helly groups CAT(0)? 
\item Is every coarsely Helly group a Helly group?
\item Is any asymptotic cone of an injective metric space injective? Note that we have a partial answer in Corollary~\ref{cor:asymptotic_cone_injective}.
\item Do Helly graphs satisfy some form of rank rigidity? Even stating a reasonable conjecture is not obvious. Here is a weak version proposed by Genevois (note that the CAT(0) version has been proved by Kent and Ricks (see~\cite[Coroallary~6.22]{kent_ricks_asymptotic_CAT0}, see also~\cite{petyt_spriano_zalloum_curtains_cat0} for another proof using curtains).
\een

\bconj[Genevois]
Let $G$ denote a group acting geometrically on a Helly graph $X$. Then one of the following holds:
\bit
\item Every asymptotic cone of $X$ has no cut-point.
\item $G$ contains a contracting isometry.
\eit
\econj

Also note that Sisto and Zalloum proved that, for any injective metric space, Morse isometries are contracting (see~\cite{sisto_zalloum_morse_contracting}). See also the survey of Zalloum about the use of injective metric spaces for hyperbolic-like behaviour of groups (\cite{zalloum_survey_curtains}).

\mk

We believe that the combinatorial dimension of a metric space is a very powerful notion, that should be more throroughly studied. Here are questions related to combinatorial dimension:

\ben
\item Is there a local characterization of combinatorial dimension? More precisely, if $X$ is a contractible metric space which has locally combinatorial dimension bounded above by some $N \in \N$, does $X$ has combinatorial dimension bounded above by $N$?
\item If $X$ is the $1$-skeleton of an $n$-dimensional CAT(0) cube complex, is the combinatorial dimension of $X$ equal to $2^{n-1}$? What is the Helly hull of $X$?
\item What is the combinatorial dimension of the standard Cayley graph of a Coxeter group?
\item What is the combinatorial dimension of a Euclidean buiding?
\item Assume that $G$ acts properly and cocompactly on a finite-dimensional injective metric space. Is $G$ Helly?
\een

\newpage

\section{Exercises} \label{sec:exercises}

\subsection{Injective spaces}

\bex Recall that a geodesic metric space is injective if and only if any family of pairwise intersecting balls has a non-empty intersection.
\ben
\item Show that the Euclidean plane $(\R^2,\ell^2)$ is not injective.
\item Show that the real line $\R$ is injective.
\item Show that the plane $(\R^2,\ell^\infty)$ is injective.
\item Can you find another norm on $\R^2$ which is injective?
\een
\eex

\bex
Let $(X,d_X),(Y,d_Y)$ denote injective metric spaces. Show that the $\ell^\infty$ product $(X \times Y,\max(d_X,d_Y))$ is injective.
\eex

\bex Recall that a metric space $X$ is injective if and only if, for any metric spaces $A \subset B$, any $1$-Lipschitz map $f:A \ra X$ has a $1$-Lipschitz extension to $B$.
\ben
\item Show that any injective metric space is geodesic.
\item Show that any injective metric space is modular, i.e. for any $x,y,z \in X$, there exists $m \in X$ such that
$$d(x,y)=d(x,m)+d(m,y), d(y,z)=d(y,m)+d(m,z), d(z,x)=d(z,m)+d(m,x).$$
\item Show that $(\R^3,\ell^1)$ is modular but not injective.\
(Hint: consider $A=\{(1,0,0),(0,1,0),(0,0,1),(1,1,1)\}$)
\een
\eex

\bex
Let us consider the (non-CAT(0)) square complex $X$ consisting of 3 squares arranged as in the corner of a $3$-cube, endowed with the standard $\ell^\infty$ metric on each square. Show that $X$ is not injective.
\eex

Recall that, if $X$ is a metric space, we may define
$$\Delta(X)=\{f:X \ra \R \st \forall x,y \in X, f(x)+f(y) \geq d(x,y)\}.$$
The injective hull of $X$ may be described as
$$E(X)=\{f \in \Delta(X) \st \forall g \in \Delta(X), g \leq f \Longrightarrow f=g\},$$
endowed with the sup metric.

\bex
Let us consider the $3$ point metric space $X=\{x_1,x_2,x_3\}$, with $d(x_i,x_j)=1$ for $i \neq j$.
\ben
\item The space $\Delta(X)$ is isometric to the subspace
$$\Delta(X)=\{x \in \R^3 \st \forall i, x_i \geq 0, \forall i \neq j, x_i+x_j \geq 1\}$$
of $(\R^3,d_\infty)$.
\item The injective hull $E(X)$ of $X$ is isometric to a tripod with endpoints $\{x_1,x_2,x_3\}$ and edge lengths $\f12$.
\item Draw $E(X)$ as the minimal subset of $\Delta(X)$.
\een
\eex

\bex
Let us consider two isometric copies $P,P'$ of the plane $(\R^2,\ell^\infty)$ glued along the half-plane
$$H=\{(x,y) \in \R^2 \st y \leq ax\},$$
for some $a \in \R$. We will prove that the gluing $X=P \cup_H P'$ is injective if and only if $a \in \{-1,0,1\}$:
\ben
\item Show that, if $a \not\in \{-1,0,1\}$, then $X$ is not injective.
\item Show that, if $a=0$, the metric space $X$ is isometric to an $\ell^\infty$ product of injective spaces. Deduce that $X$ is injective.
\item Show that, if $a=\pm 1$, the metric space $X$ is injective. Here are options:
\ben
\item Suggestion: show that $X$ can be realized as a poset, with a cofinal action by $\R$. Use then Theorem~\ref{thm:lattice_injective} to prove that $X$ is injective.
\item Other suggestion: show that $X$ can be realized as the geometric realization of a simplicial complex, with the $\ell^\infty$ orthoscheme metric, and use Theorem~\ref{thm:ell_infty_realization_poset_injective}.
\item Yet another suggestion: show that $X$ is the asymptotic cone of a Helly graph, and use Corollary~\ref{cor:asymptotic_cone_injective} to conclude.
\een
\een
\eex

\subsection{Helly graphs}

Recall that a simplicial graph $\Gamma$ is Helly if and only if it is connected, and any pairwise intersecting balls have a non-empty intersection.

\bex
Consider the graph $\Gamma$ with vertex set $\Z^2$, with an edge between $v$ and $w$ if $|v_1-w_1| \leq 1$ and $|v_2-w_2| \leq 1$. Show that $\Gamma$ is a Helly graph.
\eex

\bex
Show that the following graphs are not Helly:
\bit
\item The standard square tiling of $\R^2$.
\item The standard equilateral triangle tiling of $\R^2$.
\eit
\eex

\bex If $\Gamma_1,\Gamma_2$ are simplicial graphs, their $\ell^\infty$ product $\Gamma$ is the simplicial graph with vertex set
$$\Gamma^{(0)} = \Gamma_1^{(0)} \times \Gamma_2^{(0)},$$
with an edge between $(v_1,v_2)$ and $(w_1,w_2)$ if, for $i \in \{1,2\}$, the vertices $v_i$ and $w_i$ are equal or adjacent. Show that if $\Gamma_1$ and $\Gamma_2$ are Helly graphs, then $\Gamma$ is a Helly graph.
\eex

Recall that the Helly hull of a graph $X$ is the unique minimal Helly graph containing isometrically $X$.

\bex\

\bit
\item What is the Helly hull of a $4$-cycle?
\item What is the Helly hull of a $5$-cycle?
\item What is the Helly hull of a $6$-cycle?
\eit
\eex

\bex
Let $X$ denote a Helly graph, and let
$$P_X = \{\mbox{cliques which are intersections of balls of $X$}\} \subset {\cal P}(X),$$
with the inclusion partial order. Show that $P_X$ satisfies the assumptions of Theorem~\ref{thm:ell_infty_realization_poset_injective}, i.e. it is locally a flag semilattice.
\eex

Recall that if $X$ is a Helly graph, its injective hull (of its vertex set) can be described as the geometric realization of the poset $P_X$, with the standard piecewise $\ell^\infty$ orthosimplex metric. Also recall that the first barycentric subdivision of $X$ is the graph $X'$ with vertex set $P_X$, with an edge between points at distance $\f12$ in $E(X)=|P_X|$.

\bex
Let us consider a triangle denoted $X$.
\ben
\item What is the Helly hull $H(X)$ of $X$?
\item What is the injective hull $E(X^{(0)})$ of the vertex set of $X$?
\item What is the first barycentric subdivision $X'$ of $X$?
\een
Note how $X$ and $X'$ are different graphs.
\eex

\bex
For each of the following graphs, compute the Helly hull, and compute the injective hull of the vertex set.
\ben
\item A square.
\item A square with a diagonal.
\item Two squares glued along an edge.
\item Two complete graphs over $4$ vertices glued along an edge.
\een
\eex

\subsection{Miscellaneous}

\bex
Let us consider the union $X$ of a square and a $3$-cube along an edge (see Figure~\ref{fig:geodesics_cube_lp}). Let us endow $X$ with the piecewise $\ell^p$ metric, for $p \in [1,\infty]$. Show that, for $p \in (1,\infty)$, the unique geodesic between the opposite vertices $a$ and $b$ intersects the edge $e$ in a different point.
\eex

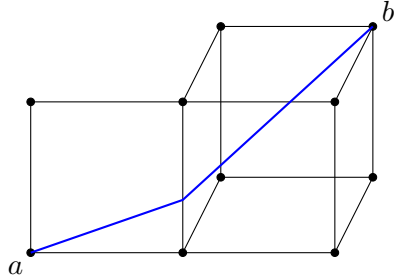
\begin{figure}[H]
\begin{center}
\begin{tikzpicture}
\def \p {0.05}
\def \t {(0.5,1)}
\def \op {0.5}
\def \gris {black!10}
\draw[fill] (0,0) circle (\p) node(a) {};
\draw[fill] (2,0) circle (\p) node(b) {};
\draw[fill] (4,0) circle (\p) node(c) {};
\draw[fill] (4,2) circle (\p) node(d) {};
\draw[fill] (2,2) circle (\p) node(e) {};
\draw[fill] (0,2) circle (\p) node(f) {};
\draw[fill] (2,0)+\t circle (\p) node(b') {};
\draw[fill] (4,0)+\t circle (\p) node(c') {};
\draw[fill] (4,2)+\t circle (\p) node(d') {};
\draw[fill] (2,2)+\t circle (\p) node(e') {};

\draw (a.center) -- (b.center) -- (c.center) -- (d.center) -- (e.center) -- (f.center) -- (a.center);
\draw (b.center) -- (e.center);
\draw (b'.center) -- (c'.center) -- (d'.center) -- (e'.center) -- (b'.center);
\draw (b.center) -- (b'.center);
\draw (c.center) -- (c'.center);
\draw (d.center) -- (d'.center);
\draw (e.center) -- (e'.center);
\node at (-0.2,-0.2) {\bfseries $a$};
\node at (4.7,3.2) {\bfseries $b$};
\draw[thick,blue] (a.center) -- (2,0.7) -- (d'.center);

\end{tikzpicture}
\end{center}
\caption{Geodesics in a cube complex with the $\ell^p$ metric}
\label{fig:geodesics_cube_lp}
\end{figure}

\bex
A metric space is $0$-hyperbolic if and only if it has combinatorial dimension $1$.
\eex

\bex\
\ben
\item Prove that any Gromov-hyperbolic graph has stable intervals.
\item Prove that any Helly graph has stable intervals, with constant $\beta=1$.
\item Prove that any median graph (i.e. $1$-skeleton of a CAT(0) cube complex) has stable intervals, with constant $\beta=1$.
\een
\eex

\addcontentsline{toc}{section}{References}
\bibliographystyle{alpha}
\bibliography{../../../bibli}

\end{document}